\documentclass[reqno]{amsart}

\usepackage{bm}
\usepackage{bm,amssymb,amsmath,enumerate,color,txfonts,mathrsfs,tipa,multirow,epsfig,indentfirst,subfig}%
\usepackage{amsmath,amsthm,amssymb,amsfonts,mathrsfs,amscd}
\usepackage{multirow}
\usepackage{graphicx}
\makeatletter
\newcommand\figcaption{\def\@captype{figure}\caption}
\newcommand\tabcaption{\def\@captype{table}\caption}
\makeatother

\theoremstyle{definition}

\newtheorem{remark}{Remark}[section]

\numberwithin{equation}{section}

\def \be{\begin{equation}}
\def \en{\end{equation}}
\def\beq{\begin{eqnarray}}
\def\eq{\end{eqnarray}}
\def\beqx{\begin{eqnarray*}}
\def\eqx{\end{eqnarray*}}
\def\l{\label}

\def\0{{\bf 0}}

\def\x{{\bf x}}

\def\12{{1\over 2}}

\def\n{{\bf n}}

\begin{document}
\title{Novel multilevel preconditioners for the systems arising from plane wave discretization of Helmholtz equations with large wave numbers}

\author{Qiya Hu}
\author{Xuan Li}

\thanks{LSEC, Academy of Mathematics and Systems Science,
Chinese Academy of Sciences, Beijing 100190, China.
(hqy@lsec.cc.ac.cn \and lixuan@lsec.cc.ac.cn). This research was supported by the Natural Science Foundation of China
G11571352.}

\maketitle

 \begin{abstract}
In this paper we are concerned with fast algorithms for the systems arising from the plane wave discretizations for two-dimensional Helmholtz equations with large wave numbers. We consider the plane wave weighted least squares (PWLS) method and the plane wave discontinuous Galerkin (PWDG) method. The main goal of this paper is to construct multilevel parallel preconditioners for solving the resulting Helmholtz systems. To this end, we first build a multilevel overlapping space decomposition for the plane wave discretization space based on a multilevel overlapping domain decomposition method. Then, corresponding to the space decomposition, we construct an additive multilevel preconditioner for the underlying Helmholtz systems. Further, we design both additive and multiplicative multilevel preconditioners with
smoothers, which are different from the standard multigrid preconditioners. We apply the proposed multilevel preconditioners with a {\it constant} coarsest mesh size to solve two dimensional Helmholtz
systems generated by PWLS method or PWDG method, and we find that the new preconditioners possess nearly stable convergence, i.e., the iteration counts of the preconditioned iterative methods (PCG or PGMRES) with the preconditioners increase very slowly when the wave number increases (and the fine mesh size decreases).
 \end{abstract}



\pagestyle{myheadings} \thispagestyle{plain} \markboth{Qiya Hu and Xuan Li}{Novel multilevel preconditioners for the systems arising from plane wave discretization of Helmholtz equations}

\section{introduction}

The plane wave method, which falls into the class of Trefftz methods \cite{ref10}, differs from the traditional finite element
method and the boundary element method in the sense that
the basis functions are chosen as exact solutions of the governing
differential equation without boundary conditions. This type of
numerical method was first introduced to solve Helmholtz equations. Examples of this
approach include the Ultra Weak Variational Formulation (UWVF) (see \cite{ref11,ref11_3}), the weighted plane wave least-squares (PWLS) method (see \cite{refhy, ref12}), the plane wave discontinuous
Galerkin methods (PWDG) (see \cite{ref20,ref21}), the plane wave Lagrangian multiplier (PWLM) method \cite{ref33, pa1} and the Variational Theory of Complex
Rays (VTCR) introduced in \cite{ref14,ref15_2, ref17}.
This kind of method can generate higher accuracy approximations than the other methods for Helmholtz equations with large wave numbers.
The plane wave discretization methods have been extended to
discretization of Maxwell's equations recently (see \cite{HMP2013, HuY2014, HMM2007}). The PWLS method has an
advantage over the other plane wave methods: the stiffness matrix associated with the PWLS method is Hermitian positive definite, so the resulting system can be solved by the PCG method.
Like the other discretization methods, the Helmholtz systems arising from the plane wave discretization are also highly ill-conditioned when the wave number is large. Comparing with many works on the plane wave discretizations, there are only a few articles (refer to \cite{ref33, refhy, pa1}) to study fast solver for the resulting Helmholtz systems.

It is well known that multilevel methods are powerful algorithms for solving the systems generated by finite element discretization of elliptic-type partial differential equations (see, for example, \cite{BankDupYes1988,BramPasXu1990,Brand1977,DrySarWid1996,Hack1985,Hiptm1998}). However, the standard multilevel methods (and domain decomposition methods) are ineffective for Helmholtz equations (and time-harmonic Maxwell's equations) with large wave numbers, unless the sizes of coarse meshes are chosen as $O(1/\omega)$ (see, for example, \cite{BramPasXu1988,refcai,Elam2001,Erl2006, Kim2002,re23, Gop2003, re1, re24, Yser1986}), where $\omega$ denotes the fixed wave number. It is clear that the restriction on the coarse mesh sizes is limiting in applications.
How to construct an effective parallel preconditioner for Helmholtz equations (and time-harmonic Maxwell's equations) with large wave numbers seems an open problem. The wave-ray multigrid method for Helmholtz equations
was proposed in \cite{Brand1997,LivBrand2006}  (a further development of this method was made in \cite{LeeMccor2000}), in which the approximations of oscillatory error components were transformed into the approximations of
smooth ray envelope functions by using the exponential interpolations. The wave-ray multigrid method can improve the performance of the standard multigrid methods for Helmholtz equations with large wave numbers. Recently, a kind of successive preconditioner based on a decomposition of the domain into strips was proposed in \cite{stddm, sweeppml} to solve Helmholtz equations with large wave numbers. The preconditioners can be viewed as physically-based
approximations of direct solvers. It has been shown that such kind of preconditioner possesses the optimal convergence independent of the mesh sizes \cite{stddm}, which is a very important result in the solution method for Helmholtz equations with large wave numbers.

In the present paper, we consider the PWLS method and the PWDG method for the discretization of Helmholtz equations
in two dimensions, and explore a new way to construct multilevel preconditioners for the resulting Helmholtz systems. At first we design a multilevel overlapping domain decomposition method to build a multilevel space decomposition for the plane wave discretization space. Then, based on the space decomposition, we construct an additive multilevel overlapping preconditioner for the underlying Helmholtz systems. Finally, we replace the solvers in the previous preconditioner by block Jacobi-type smoothers to get cheaper (both additive and multiplicative) multilevel overlapping preconditioners. The multilevel overlapping preconditioners with smoothers are different from the standard multigrid preconditioners, since the space decomposition defining such new preconditioners has different overlapping structure from the one corresponding to the standard multigrid preconditioners. We apply the proposed preconditioners to solve Helmholtz
systems generated by PWLS method or PWDG method. Numerical results indicate that the new preconditioners possess nearly stable convergence, i.e., the iteration counts of the corresponding iterative methods (PCG or PGMRES) increase very slowly when the wave number increases (and the mesh size decreases), without the limiting condition mentioned in the last paragraph. In particular, the multilevel overlapping preconditioners with smoothers possess almost optimal convergence.

The paper is organized as follows: In Section 2, we recall the PWLS method and the PWDG method for Helmholtz equations. In section 3, we design a multilevel space decomposition of the solution space and describe the corresponding additive multilevel preconditioner. An additive multilevel overlapping preconditioner with smoothers
is introduced in Section 4. In Section 5, we define several multiplicative variants of the additive multilevel overlapping preconditioner with smoothers. In Section 6, we apply the proposed preconditioners to solve several Helmholtz systems and report some numerical results.

\section{Plane wave methods for Helmholtz equations}

For convenience, we only consider the two-dimensional case in this paper.
In this section, we briefly review the plane wave methods for Helmholtz equations. At first the original problem
to be solved is defined. Then the variational formulations are given out in detail.
\subsection{The reference problem}
Firstly, we present the mathematical model of Helmholtz equations. Let $\Omega$ be a
bounded and connected Lipschitz domain in two dimensions. We consider Helmholtz equations with Robin boundary conditions.
\begin{eqnarray}
\left\{\begin{array}{ll} -\Delta u-\omega^2u=0& \text{in}\quad
\Omega,\\
(\partial_\text{\bf n}+i\omega)u=g& \text{on}\quad\gamma=\partial\Omega,
\end{array}\right.
\label{eq11}
\end{eqnarray}
where $\partial_{\bf n}$ and $\omega$ denote the outer normal derivative and the angular frequency.

Let $\Omega$ be divided into a partition as follows:
\[
\overline{\Omega}=\bigcup_{k=1}^N\overline{E}_k,\quad
E_k\cap E_j=\emptyset
 \quad\text{  for  }k\not=j.
 \]
We assume that the subdomains $E_1,E_2,\cdots,E_{N}$ are geometrical conforming, i.e., the intersection of any two adjoining subdomains
is just the common vertex or the common edge of them. Here, we do not require that the intersection of two adjoining elements
is a straight line segment. In practice, the partition is a mesh of domain, and $E_1,\cdots,E_N$ are the elements. As usual, we assume that $\{E_k\}$ is quasi-uniform and regular.
Let $\mathcal{T}_h$ denote the set of the elements $E_1,\cdots,E_N$,
where $h$ is the size of the elements. Define
$$ \gamma_{kj}=\partial E_k\cap\partial E_j \quad\quad(\mbox{when}~~ E_k~~ \mbox{and}~~ E_j~~ \mbox{are}~~ \mbox{adjoining}) $$
 and
$$ \gamma_k=\partial E_k\cap\partial\Omega \quad\quad(\mbox{if}~~ E_k~~ \mbox{closes}~~ \partial\Omega).$$

Let $V(E_k)$ denote the space of the functions which verify
Helmholtz's homogeneous equation (\ref{eq11}) on the element
$E_k$:
\begin{equation}
V(E_k)=\{v_k\in H^1(E_k);~ \Delta v_k+\omega^2 v_k=0\}.
\end{equation}
Define
$$V({\mathcal T}_h)=\prod\limits_{k=1}^NV(E_k),$$
with the natural scalar product
$$ (u,v)_V =\sum\limits_{k=1}^N\int_{E_k}u_k\cdot\overline{v}_k \ d{\bf x},
~~\forall u,v \in V({\mathcal T}_h).
$$

\subsection{The PWLS method}
In this subsection, we review the PWLS method introduced in \cite{ref12} and \cite{refhy}.

Set $u|_{E_k}=u_k$ $(k=1,\cdots,N$). Then
the reference problem to be solved consists in finding the local
acoustic pressures $u_k\in H^1(E_k)$ such that
\begin{equation}
\left\{\l{eq1}
\begin{array}{rrll}
-\Delta u_k-\omega^2u_k&=&0 \quad&\text{in}\quad E_k,  \\
(\partial_{\n}+i\omega)u&=&g\quad&\text{on}\quad\gamma_k~~(if~~\gamma_k\not=\emptyset),
\end{array}
\right.
\end{equation}
and
\begin{eqnarray}
\left\{\begin{array}{rrll}u_k-u_j&=&0 & \text{over}\quad \gamma_{kj},\\
\partial_{\text{\bf n}_k}u_k+\partial_{\text{\bf n}_j}u_j&=&0& \text{over}\quad
\gamma_{kj}
\end{array}\right.\quad\quad(k\neq j;~ k,j=1,2,\cdots,N).
\label{eq3}
\end{eqnarray}

Let $\alpha$ and $\beta$ be two given positive real numbers to be specified later. Corresponding to the boundary condition in (\ref{eq1}) and the
interface continuity condition (\ref{eq3}), we define the functional
\begin{equation}
\begin{split}
\quad J(v)&=\sum_{k=1}^N\int_{\gamma_k}|(\partial_{\bf
n}+i\omega )v_k-g|^2ds \\&  + \sum_{j\not=k}
\bigg(\alpha\int_{\gamma_{kj}}
|v_k-v_j|^2ds+\beta\int_{\gamma_{kj}}|\partial_{\text{\bf n}_k}
v_k+\partial_{\text{\bf n}_j}v_j|^2ds\bigg),~~v\in V({\mathcal
T}_h). \label{minJ}
\end{split}
\end{equation}

It is clear that $J(v)\geq 0$. Consider the minimization problem:
find $u\in V({\mathcal T}_h)$ such that
\begin{equation}
J(u)=\min\limits_{v\in V({\mathcal T}_h)}J(v)\label{min}
\end{equation}
If $u$ is the solution of the problem (\ref{eq11}),
i.e., $u\in V({\mathcal T}_h)$ satisfies the boundary condition in
(\ref{eq1}) and the interface continuity condition (\ref{eq3}), then
we have $J(u)=0$, which implies that $u$ is also the solution of the
minimization problem (\ref{min}).

Define the sesquilinear form $a(\cdot,\cdot)$ by
\begin{eqnarray}
 &&a(u,v)=\sum_{k=1}^N\int _{\gamma_k}((\partial_{\bf
n}+i\omega )u_k)\cdot\overline{(\partial_{\bf n}+i\omega) v_k}ds \cr
&&\quad\quad\quad \quad +
 \sum_{j\not=k} \bigg( \alpha\int_{\gamma_{kj}}
(u_k-u_j)\cdot\overline{(v_k-v_j)}ds \cr&&\quad \quad
+\beta\int_{\gamma_{kj}}(\partial_{\text{\bf n}_k}
u_k+\partial_{\text{\bf n}_j} u_j
)\cdot\overline{(\partial_{\text{\bf n}_k} v_k+\partial_{\text{\bf
n}_j} v_j )} ds\bigg), ~~\forall v\in V({\mathcal T}_h), \label{eq7}
\end{eqnarray}
and define the functional $\mathcal{L}(\cdot)$ by
\begin{eqnarray}
&&\mathcal{L}(v)=\sum_{k=1}^N\int
_{\gamma_k}g\cdot\overline{(\partial_{\bf n}+i\omega)
v_k}ds~~\quad\forall v\in V({\mathcal T}_h).
\end{eqnarray}

The variational problem associated with the minimization problem
(\ref{min}) can be expressed as:

\begin{eqnarray}
\left\{\begin{array}{ll} Find \ u\in V({\mathcal T}_h), s.t. \\
a(u,v)=\mathcal{L}(v), \quad \forall v\in V({\mathcal T}_h).
\end{array}\right.
\label{eq30}
\end{eqnarray}

The reference problem (\ref{eq1}) and (\ref{eq3}) is equivalent to the
new variational problem (\ref{eq30}) (see \cite{refhy} Theorem 3.1). In applications, we usually choose the two parameters
in (\ref{minJ}) as $\alpha=\omega^2$ and $\beta=1$.

\subsection{The PWDG method}
In this subsection, we review the PWDG method introduced in \cite{ref21}.

Let $u$ and ${\bf \sigma}$ be a piecewise smooth function and vector field on ${\mathcal T}_h$ respectively. On $\gamma_{kj}$, we define
\begin{eqnarray*}
\text{the averages:} & \{u\}=\frac{1}{2}(u_k+u_j),\quad \{{\bf\sigma}\}=\frac{1}{2}({\bf \sigma}_k+{\bf \sigma}_j),\\
\text{the jumps:} & \quad[u]=u_k{\bf n}_k+u_j{\bf n}_j,\quad [{\bf \sigma}]={\bf \sigma}\cdot{\bf n}_k+{\bf \sigma}\cdot{\bf n}_j.
\end{eqnarray*}
Set
$$ \mathcal{F}_h^I=\bigcup_{k\neq j}\gamma_{kj}\quad\mbox{and}\quad\mathcal{F}_h^B=\bigcup^N_{k=1}\gamma_k.$$

With these definitions, we can write the PWDG method as follows:
\begin{eqnarray}
\left\{\begin{array}{ll} Find \ u\in V({\mathcal T}_h), \ s.t. \\
a(u,v)=\mathcal{L}(v),~~\forall v\in V({\mathcal T}_h),
\end{array}\right.\label{PWDG}
\end{eqnarray}
where (see \cite{ref21})
\begin{equation}
\begin{split}
a(u,v)&=\int_{\mathcal{F}_h^I}\bigg(\{u\}[\overline{\nabla v}]-\frac{\beta}{i\omega}[\nabla u][\overline{\nabla v}]-\{\nabla u\}[\overline{v}]+\alpha\cdot i\omega[u][\overline{v}]\bigg)ds\\
&+\int_{\mathcal{F}_h^B}\bigg((1-\delta)u\overline{\nabla v\cdot {\bf n}}-\frac{\delta}{i\omega}\nabla u\cdot{\bf n}\overline{\nabla v\cdot{\bf n}}-\delta\nabla u\cdot{\bf n}\overline{v}+(1-\delta)i\omega u\overline{v}\bigg)ds,
\end{split}
\end{equation}
and
\begin{equation}
\mathcal{L}(v)=
\int_{\mathcal{F}_h^B}\bigg(-\frac{\delta}{i\omega}g\overline{\nabla v\cdot{\bf n}}+(1-\delta)g\overline{v}\bigg)ds.
\end{equation}
Here $\alpha,\beta$ and $\delta$ are given positive parameters. The simplest choice of the parameters in the above two expressions is $\alpha=\beta=\delta={1\over 2}$.

\subsection{Discretization of the variational formulations}
Before building discrete variational problems, we need to approximate the space $V({\mathcal T}_h)$ by a suitable finite dimensional subspace, which is spanned by some plane wave basis functions, i.e., solutions of homogeneous Helmholtz equation without boundary condition.

For convenience, we assume that the number of plane wave basis functions equals a same positive integer $p$ for every elements $\Omega_k$. Let $y_{l}$ be the wave
shape functions, which satisfy
\begin{eqnarray}
\left\{\begin{array}{ll}  y_{l}({\bf x})=e^{i\omega(\boldsymbol{\alpha_l}\cdot{\bf x})},~~{\bf
x}\in\Omega,\\
\boldsymbol{\alpha_l}\cdot \boldsymbol{\alpha_l}=1,\\
l\neq s \rightarrow
\boldsymbol{\alpha_{l}}\neq\boldsymbol{\alpha_{s}},
\end{array}\right.
\label{eq31}
\end{eqnarray}
where $\boldsymbol{\alpha_l}~(l=1,\cdots,p)$ are unit wave
propagation directions to be specified later. The plane wave basis functions can be defined as
\begin{eqnarray}\quad\quad\phi^{(k)}_{l}({\bf x})=\left\{\begin{array}{ll} y_{l}({\bf x}),~~{\bf
x}\in E_k,\\
\quad 0,\quad {\bf
x}\notin E_k
\end{array}\right.\quad~~(k=1,\cdots,N;~l=1,\cdots,p).
\label{new}
\end{eqnarray}
Thus the space $V({\mathcal T}_h)$ is discretized by the subspace
\begin{equation}
V_p({\mathcal T}_h)=span\bigg\{\phi^{(k)}_{l}:~k=1,\cdots,N;~l=1,\cdots,p
\bigg\}. \label{eq32}
\end{equation}

During numerical simulations, the directions of the
wave vectors of these wave functions, for two-dimensional problems,
are uniformly distributed as follows: $$
\boldsymbol{\alpha_l}=\begin{pmatrix}(cos(2\pi(l-1)/p)\\sin(2\pi(l-1)/p))\end{pmatrix}~~~(l=1,\cdots,p).$$

Let $V_p({\mathcal T}_h)$ be the plane wave space defined above. Then
the discrete variational problems associated with (\ref{PWDG}) and (\ref{eq30}) can be
described as follows:
\begin{eqnarray}
\left\{\begin{array}{ll} Find \ u_h\in V_p({\mathcal T}_h), \ s.t. \\
a(u_h,v_h)=\mathcal{L}(v_h),~~\forall v_h\in V_p({\mathcal T}_h).
\end{array}\right.
\label{eq33}
\end{eqnarray}

Let $A: V_p({\mathcal T}_h)\rightarrow V_p({\mathcal T}_h)$ be the discrete operator defined by the sesquilinear form $a(\cdot,\cdot)$.  The discrete variational problem (\ref{eq33}) can be written in the operator form
\be
Au_h=f_h,~~~u_h\in V_p({\mathcal T}_h)\label{operator form}.
\en

Let ${\mathcal A}$ be the stiffness matrix generated by the
sesquilinear form $a(\cdot,\cdot)$ on the space $V_p({\mathcal
T}_h)$, and let $b$ denote the vector associated with $\mathcal{L}(v_h)$. Namely, the entries of the matrix ${\mathcal
A}$ are computed by $a_{k,j}^{l,m}=a(\phi^{(j)}_{m},\phi^{(k)}_{l})$; and the
complements of the vector $b$ are defined as
$b_{k,l}=\mathcal{L}(\phi^{(k)}_{l})$. The discretized problem
(\ref{eq33}) leads to the algebraic system below: \be {\mathcal
A}X=b, \label{eq34} \en where
$X=(x_{11},x_{12},\cdots,x_{1p},x_{21},\cdots,x_{2p},\cdots,x_{N1},\cdots,x_{Np})^t
\in \mathbb{C}^{pN}$ is the unknown vector.

In general the system (\ref{eq33}) is solved by an iterative method, for example, the preconditioned GMRES method or the PCG method. In this paper, we solve the system arising from the PWDG method by preconditioned GMRES method, and solve the system arising from the PWLS method by PCG method since the system of the PWLS method is Hermitian positive definite. Notice that implementation of an iterative step in PCG method is cheaper than that in the preconditioned GMRES method. We need
to construct an efficient preconditioner ${\mathcal B}$ for the
matrix ${\mathcal A}$, and solve the equivalent
system
\begin{equation}
{\mathcal B}^{-1}{\mathcal A}X={\mathcal B}^{-1}b.\label{system}
\end{equation}

The main goal of this paper is to construct efficient multilevel preconditioners ${\mathcal B}$, especially multilevel preconditioners with overlapping smoothers. In order to make the ideas easily understood,
we first construct a basic preconditioner directly from multilevel overlapping domain decompositions, and then we define multilevel preconditioners with overlapping smoothers based on the basic preconditioner.
For convenience, we shall describe the preconditioners in operator forms, instead of matrix forms.

\section{A preconditioner based on multilevel overlapping domain decomposition}
In this section, we construct an additive multilevel preconditioner $B$ for the operator $A$ based on overlapping domain decompositions.

\subsection{A multilevel overlapping space decomposition}
Let $N_0$ be a fixed positive integer, which is independent of $\omega$, $h$ and $p$. For simplicity of exposition, we use $D$ to denote a generic domain that
is the union of some elements in ${\mathcal T}_h$, where $D$ can be the domain $\Omega$ itself or a subdomain of $\Omega$.

Let $D$ be decomposed into the union of non-overlapping subdomains $D_1,D_2,\cdots,D_{N_0}$
such that: (1) each subdomain $D_r$ is just the union of several elements in
${\mathcal T}_h$; (2) the subdomains $D_1,D_2,\cdots,D_{N_0}$ are quasi-uniform, regular and geometrical conforming (refer to Subsection 2.1).
Here, we do not require that the intersection of two adjoining subdomains is a straight line segment. Then
$D_1,\cdots,D_{N_0}$ can be viewed as coarse elements of $D$ and they
constitute a (coarse) finite element partition ${\mathcal T}^D_d$ of $D$, where $d$ denotes the size of these elements.

Based on the partition ${\mathcal T}^D_d$, we can define an overlapping domain decomposition of $D$ as usual.
For a constant $\theta_0\in [{1\over 2},1]$,
we enlarge each coarse element $D_r$ by the thickness $\theta_0d$, and generate
a larger domain $\tilde{D}_r$ satisfying: (1) $D_r\subset\tilde{D}_r\subset D$; (2) $\tilde{D}_r$ is just the union of some (fine) elements in ${\mathcal T}_h$;
(3) the distance between the internal boundaries $\partial \tilde{D}_r\backslash\partial D$ and $\partial D_r\backslash\partial D$ is about $\theta_0d$.
Then
$$ \bar{D}=\bigcup_{r=1}^{N_0}\tilde{D}_{r} $$
constitutes an overlapping domain decomposition of $D$ with ``large overlap". For convenience, we call the parameter $\theta_0$ as ``overlapping degree".
When $\theta_0=1$ (rep. $\theta_0={1\over 2}$), each subdomain $\tilde{D}_r$ is the union of $D_r$ itself and all the neighboring coarse elements (rep.
the half of every neighboring coarse elements) with it. Thus, the case with $\theta_0=1$ (rep. $\theta_0={1\over 2}$) is called ``complete overlap" (rep.
``half overlap"). We point out that the case with a small $\theta_0$, i.e., small overlap (for example, $\theta_0={h\over d}$)
is not considered in this paper, since the numerical results for this case are not satisfactory (see Table \ref{table-new2} in Section 6).

For convenience, the above process to generate the coarse elements $\{D_r\}$ and the overlapping subdomains $\{\tilde{D}_r\}$ from $D$ is called
a ``decomposition operation" of $D$. The subdomain $\tilde{D}_r$ is called the ``enlarged subdomain" of $D_r$.

When $D$ is just $\Omega$ itself, we let ${\mathcal T}^{\Omega}_{d_0}$
denote the set of the resulting coarse elements $\Omega_1,\cdots,\Omega_{N_0}$, where $d_0$ is the size of the elements $\Omega_1,\cdots,\Omega_{N_0}$.
Moreover, we use ${\mathcal S}_0$ to denote the set of the ``enlarged subdomains" $\tilde{\Omega}_1$, $\tilde{\Omega}_2$,
$\cdots,\tilde{\Omega}_{N_0}$.

For each subdomain $D\in{\mathcal S}_0$, let ${\mathcal T}^D_{d_1}$ be the set of the coarse elements $D_1,\cdots,D_{N_0}$ defined by the ``decomposition operation" of $D$,
where $d_1$ denote the size of $D_1,\cdots,D_{N_0}$. Let
$\bar{D}=\bigcup_{r=1}^{N_0}\tilde{D}_{r}$ denote the overlapping domain decomposition of $D$, where $\tilde{D}_{r}$ is the ``enlarged subdomain" of $D_r$. With all the
``enlarged subdomains" at $1$th-level, define the set
$$ {\mathcal S}_1=\{\tilde{D}_{r}:~~r=1,\cdots,N_0;~~\mbox{for~every}~D\in {\mathcal S}_0\}.$$

We can repeat the above process. Let $J\geq 1$. For an integer $j$ satisfying $1\leq j\leq J$, we assume that the set ${\mathcal S}_{j-1}$ consisting of overlapping subdomains of $\Omega$ has been defined.
For each subdomain $D\in{\mathcal S}_{j-1}$, we use ${\mathcal T}^D_{d_j}$ to denote the set of the coarse elements $D_1,\cdots,D_{N_0}$ defined by the ``decomposition operation" of $D$, with $d_j$ being
the size of the subdomains $D_1,\cdots,D_{N_0}$. Let $\tilde{D}_{r}$ be the ``enlarged subdomain" of $D_r$, and let
$\bar{D}=\bigcup_{r=1}^{N_0}\tilde{D}_{r}$ denote the resulting overlapping domain decomposition of $D$. Define the set of $j$th-level ``enlarged subdomains" as
$$ {\mathcal S}_j=\{\tilde{D}_{r}:~~r=1,\cdots,N_0;~~\mbox{for~every}~D\in {\mathcal S}_{j-1}\}\quad\quad(j=1,\cdots,J).$$

We would like to point out that the numbers of the coarse elements generated by ``decomposition operation" of two different subdomains may be different in applications, here the choice of the same number $N_0$ of
coarse elements is only to simplify the description. When choosing $N_0$ properly, we have $d_0>d_1>\cdots>d_J>h$. Then the number of {\it fine} elements contained in each $K\in {\mathcal S}_j$ decreases rapidly when $j$ increases.

Corresponding to a ``decomposition operation" of a subdomain, we can build a local space decomposition on the subdomain.

As in Section 2, let $y_l({\bf x })$ denote the plane wave shape function $e^{i\omega ({\bf \alpha_l}\cdot{\bf x })}$ ($l=1,\cdots,p$). Let ${\mathcal Q}_p$
be the space consisting of the $p$ plane wave shape functions, i.e.,
$$ {\mathcal Q}_p=span\{y_l:~l=1,\cdots,p\}. $$
Define the {\it coarsest} plane wave space on $\Omega$ as
$$ V_p({\mathcal T}^{\Omega}_{d_0})=\{v\in L^2(\Omega):~~v|_K\in {\mathcal Q}_p~~\mbox{for every}~~ K\in {\mathcal T}^{\Omega}_{d_0}\}. $$
Similarly, for each $D\in {\mathcal S}_{j-1}$ with $j\geq 1$, define the {\it coarse} plane wave space on $D$ by
$$ V_p({\mathcal T}^D_{d_j})=\{v\in L^2(\Omega):~supp~~ v\subset D;~~v|_K\in {\mathcal Q}_p~~\mbox{for every}~~ K\in {\mathcal T}^D_{d_j}\}\quad(j=1,\cdots,J), $$
namely, $V_p({\mathcal T}^D_{d_j})$ is the plane wave finite element space associated with the coarse partition ${\mathcal T}^D_{d_j}$. Notice that the spaces
$V_p({\mathcal T}^{\Omega}_{d_0})$ and $ {V}_p({\mathcal T}^D_{d_j})$ ($j=1,\cdots, J$) have the dimension $N_0p$ and possess the same structure with the original plane
wave finite element space $V_p({\mathcal T}_h)$ defined in Subsection 2.4.

For a subdomain $K$ that is the union of some fine elements in ${\mathcal T}_h$, we always use ${\mathcal T}^K_h$ to denote the restriction of the original partition ${\mathcal T}_h$ on $K$,
and define the {\it fine} plane wave space on $K$ by
$$ V_p({\mathcal T}^K_h)=\{v\in V_p({\mathcal T}_h):~~supp~~v\subset K\}.$$

As in the standard overlapping domain decomposition method, we can obtain the initial space decomposition on $\Omega$ (here we can easily define weight functions satisfying the partition of unity,
since we do not require the continuity of functions in the considered spaces)
\begin{equation}
 V_p({\mathcal T}_h)=V_p({\mathcal T}^{\Omega}_{d_0})+\sum_{r=1}^{N_0}V_p({\mathcal T}^{\tilde{\Omega}_r}_h)=V_p({\mathcal T}^{\Omega}_{d_0})+\sum_{D\in {\mathcal S}_0}V_p({\mathcal T}^{D}_h). \label{3.intial_decomposition}
\end{equation}
Similarly, for each $D\in {\mathcal S}_{j-1}$ with $j\geq 1$, we can build the local space decomposition on $D$
\begin{equation}
 V_p({\mathcal T}^D_h)=V_p({\mathcal T}^D_{d_j})+\sum_{r=1}^{N_0}V_p({\mathcal T}^{\tilde{D}_r}_h).\label{3.local_decomposition}
\end{equation}

Set $j=1$ in (\ref{3.local_decomposition}), and substituting the resulting decomposition into (\ref{3.intial_decomposition}), yields
\beqx
 V_p({\mathcal T}_h)&=&V_p({\mathcal T}^{\Omega}_{d_0})+\sum_{D\in {\mathcal S}_0}V_p({\mathcal T}^D_{d_1})+\sum_{D\in {\mathcal S}_0}\sum_{r=1}^{N_0}V_p({\mathcal T}^{\tilde{D}_r}_h)\\
&=&V_p({\mathcal T}^{\Omega}_{d_0})+\sum_{D\in {\mathcal S}_0}V_p({\mathcal T}^D_{d_1})+\sum_{D\in {\mathcal S}_1}V_p({\mathcal T}^D_h). \eqx
Combining the above decomposition with (\ref{3.local_decomposition}) for $j=2,\cdots,J$, and using the relation
$$ \sum_{D\in {\mathcal S}_{j-1}}\sum_{r=1}^{N_0}V_p({\mathcal T}^{\tilde{D}_r}_h)=\sum_{D\in {\mathcal S}_j}V_p({\mathcal T}^D_h)\quad(j\geq 2), $$
we recursively obtain the multilevel space decomposition
\begin{equation}
V_p({\mathcal T}_h)=V_p({\mathcal
T}^{\Omega}_{d_0})+\sum_{j=1}^{J}\sum_{D\in {\mathcal S}_{j-1}}V_p({\mathcal
T}^D_{d_j})+\sum_{K\in {\mathcal S}_{J}}V_p({\mathcal T}^{K}_h)
.\label{decom-new}
\end{equation}

For ease of notation, we would like to give a terser expression of the above space decomposition.

For convenience, we write $V_p({\mathcal T}^{\Omega}_{d_0})$ as $V_p({\mathcal T}_{d_0})$. For $j=1,\cdots,J$, define the {\it set} of $j$th-level {\it coarse} elements
$$ {\mathcal T}_{d_j}=\bigcup_{D\in {\mathcal S}_{j-1}}{\mathcal T}^D_{d_j} $$
and $j$th-level {\it coarse} space
$$ V_p({\mathcal T}_{d_j})=\sum\limits_{D\in {\mathcal S}_{j-1}}V_p({\mathcal T}^{D}_{d_j}). $$
Notice that, for $j\geq 1$, the set ${\mathcal T}_{d_j}$ does not constitute a (coarse) finite element partition of $\Omega$ since the elements in
${\mathcal T}^{D}_{d_j}$ may be overlapping with the elements in ${\mathcal T}^{D'}_{d_j}$ when $D$ is different from $D'$.

Moreover, we define the set of $J$th-level {\it fine} elements
$$ \tilde{\mathcal T}^J_h=\bigcup_{K\in {\mathcal S}_{J}}{\mathcal T}^{K}_h $$
and $J$th-level {\it fine} space
$$ V_p(\tilde{\mathcal T}^J_h)=\sum_{K\in {\mathcal S}_{J}}V_p({\mathcal T}^{K}_h). $$
Also, the set $\tilde{\mathcal T}^J_h$ is not a (fine) finite element partition of $\Omega$.

Therefore, the space decomposition (\ref{decom-new}) can be simplified as
\begin{equation}
V_p({\mathcal T}_h)=V_p(\tilde{\mathcal T}^J_h)
+\sum_{j=0}^{J}V_p({\mathcal T}_{d_j}).\label{decom6}
\end{equation}

In the rest of this paper, we construct several multilevel preconditioners for $A$ based on the above multilevel space decomposition.

\subsection{A multilevel overlapping preconditioner}

In this subsection, we construct a basic preconditioner of $A$ by the multilevel space decomposition (\ref{decom6}).

Let $A_0: V_p({\mathcal T}_{d_0})\to V_p({\mathcal T}_{d_0})$ be the restriction of the discrete operator $A$ on the coarsest space $V_p({\mathcal T}_{d_0})$, namely,
\[
(A_0v_0,w_0)=a(v_0,w_0),\quad v_0\in V_p(\mathcal{T}_{d_0}),\quad \forall w_0\in V_p({\mathcal T}_{d_0}).
\]
As usual, $A_0$ is called the {\it coarsest} solver.

Let $j=1,\cdots,J$. For $D\in {\mathcal S}_{j-1}$, let $V_p({\mathcal T}^{D}_{d_j})$ be the local coarse spaces defined in the last subsection. Define $j$th-level local coarse
solvers $A^{D}_{d_j}: V_p({\mathcal T}^{D}_{d_j})\rightarrow V_p({\mathcal T}^{D}_{d_j})$ by
$$
(A^{D}_{d_j}v,w)=a(v,w),\quad v\in V_p({\mathcal T}^{D}_{d_j}),~\forall
w\in V_p({\mathcal T}^{D}_{d_j})\quad\quad (1\leq j\leq J;~D\in {\mathcal S}_{j-1}). $$
Then we define {\it inexact} solver $B_j:~V_p({\mathcal T}_{d_j})\to V_p({\mathcal T}_{d_j})$ at $j$th-level coarse space as:
\[
B^{-1}_j=\sum_{D\in {\mathcal S}_{j-1}}(A^{D}_{d_j})^{-1}Q^{D}_{d_j}\quad\quad (j=1,\cdots,J),
\]
where $Q^{D}_{d_j}: V_p({\mathcal T}_{d_j})\rightarrow V_p({\mathcal T}^{D}_{d_j})$ denote the $L^2$ projectors. Notice that
the operator $B_j$ can be viewed as a ``block-diagonal" preconditioner for the restriction of $A$ on $j$th-level coarse subspace
$V_p(\mathcal{T}_{d_j})$, where the order of each ``block" equals $N_0p$.

Similarly, for each $K\in {\mathcal S}_J$, define $J$th-level {\it local} solver $\tilde{A}^{K}_{J}: V_p(\tilde{\mathcal T}^{K}_h)
\rightarrow V_p(\tilde{\mathcal T}^{K}_h)$ by
$$ (\tilde{A}^{K}_{J}v,w)=a(v,w),\quad v\in V_p(\tilde{\mathcal T}^{K}_h),~\forall
w\in V_p(\tilde{\mathcal T}^{K}_h)\quad\quad(K\in {\mathcal S}_J),$$
and define $J$th-level fine {\it inexact} solver $\tilde{B}_J:~V_p(\tilde{\mathcal T}^J_h)\rightarrow V_p(\tilde{\mathcal T}^J_h)$ as
$$ \tilde{B}^{-1}_J=\sum_{K\in {\mathcal S}_J}(\tilde{A}^{K}_J)^{-1}\tilde{Q}^{K}_{J},$$
where $\tilde{Q}^{K}_{J}: V_p(\tilde{\mathcal T}^J_h)\rightarrow V_p(\tilde{\mathcal T}^{K}_h)$ denote the $L^2$ projectors. It is clear that
$\tilde{B}_J$ is also a ``block-diagonal" preconditioner for the restriction of $A$ on the fine subspace $V_p(\tilde{\mathcal T}^J_h)$.

Finally, corresponding to the multilevel space decomposition (\ref{decom6}), an additive multilevel preconditioner $B:V_p(\mathcal{T}_h)\to V_p(\mathcal{T}_h)$ is naturally defined as
\begin{equation}
B^{-1}=A^{-1}_0Q_0+\sum_{j=1}^{J}B^{-1}_jQ_j+\tilde{B}^{-1}_J\tilde{Q}_J, \label{3.preconditioner}
\end{equation}
where $Q_j$ ($j=0,\cdots,J$) and $\tilde{Q}_J$ denote the $L_2$ projectors into $V_p(\mathcal{T}_{d_j})$ and $V_P(\tilde{\mathcal T}^J_h)$, respectively.

The action of $B^{-1}$ can be described by the following algorithm.
\vskip 0.2in

{\bf Algorithm 3.1}. For $\xi\in V_p(\mathcal{T}_h)$, the function $u_{\xi}=B^{-1}\xi\in V_p(\mathcal{T}_h)$ can be obtained as follows:

Step 1. Computing $u_0\in V_p(\mathcal{T}_{d_0})$ by
\[
(A_0u_0,v_0)=(\xi,v_0),\quad \forall v_0\in V_p(\mathcal{T}_{d_0});
\]

Step 2. For $j=1,\cdots, J$, computing $u_{d_j}\in V_p({\mathcal T}_{d_j})$ in parallel by
\[
(B_ju_{d_j},v)=(\xi,v),\quad \forall v\in V_p({\mathcal T}_{d_j});
\]

Step 3. Computing $\tilde{u}_h^J\in V_P(\tilde{\mathcal T}^J_h)$ by
\[
(\tilde{B}_J\tilde{u}^J_h,v_h)=(\xi,v_h),\quad \forall v_h\in V_P(\tilde{\mathcal T}^J_h);
\]

Set
\[
u_\xi=u_0+\sum_{j=1}^{J}u_{d_j}+\tilde{u}^J_h.
\]

By the definitions of the solvers $B_j$ ($j=1,\cdots,J$) and $\tilde{B}_J$, Step 2-Step 3 in {\bf Algorithm 3.1} can be implemented in smaller spaces
(otherwise, {\bf Algorithm 3.1} has no significance).

The action of $B^{-1}_j$ ($j=1,\cdots,J$) appeared in Step 2 of {\bf Algorithm 3.1} can be described by the following algorithm

\vskip 0.2in

{\bf Algorithm 3.2}. For $\eta\in V_p({\mathcal T}_{d_j})$, the function $w_{\eta}=B^{-1}_j\eta\in V_p({\mathcal T}_{d_j})$ can be obtained
by two steps:

Step 1. For $D\in {\mathcal S}_{j-1}$, computing $w^{D}_{d_j}\in V_p({\mathcal T}^{D}_{d_j})$ in parallel by
$$ a(w^{D}_{d_j}, v)=(\eta,v),\quad\forall v\in V_p({\mathcal T}^{D}_{d_j});$$

Step 2. Set
$$ w_{\eta}=\sum_{D\in {\mathcal S}_{j-1}}w^{D}_{d_j}. $$

Similarly, the action of $\tilde{B}^{-1}_J$ appeared in Step 3 of {\bf Algorithm 3.1} can be described by the following algorithm

\vskip 0.2in

{\bf Algorithm 3.3}. For $\eta\in V_p(\tilde{\mathcal T}^J_h)$, the function $\tilde{w}_{\eta}=\tilde{B}^{-1}_J\eta\in V_p(\tilde{\mathcal T}^J_h)$ can be obtained by two steps:

Step 1. For $K\in {\mathcal S}_J$, computing $\tilde{w}^{K}_h\in V_p(\tilde{\mathcal T}^{K}_{h})$ in parallel by
$$ a(\tilde{w}^{K}_h, v_h)=(\eta,v),\quad\forall v_h\in V_p(\tilde{\mathcal T}^{K}_h);$$

Step 2. Set
$$ \tilde{w}_{\eta}=\sum_{K\in {\mathcal S}_J}\tilde{w}^{K}_h. $$

\vskip 0.1in
In applications, the action of $B^{-1}$ is implemented in parallel by Step 1 in {\bf Algorithm 3.1}-{\bf Algorithm 3.3}.
\begin{remark} Notice that the dimension of the coarsest space $V_p({\mathcal T}_{d_0})$ and each local ``coarse" space $V_p({\mathcal T}^{D}_{d_j})$ equals $N_0p$. Moreover, the number of fine elements contained in
$K\in {\mathcal S}_J$ monotonically decreases when $J$ increases (assume that $N_0$ is chosen in a suitable rule). Therefore, in order to
guarantee that every local space has almost the same dimension, we should choose $J$ to be large enough such
that each domain $K\in {\mathcal S}_J$ contains almost $N_0$ fine elements in ${\mathcal T}_h$.
 Then each subproblem needed to be solved in Step 1 of {\bf Algorithm 3.2}-{\bf Algorithm 3.3} has nearly $N_0p$ unknowns only.
\end{remark}

\subsection{Further discussions on the proposed multilevel method}

In this subsection we first give some comparisons between the proposed multilevel method and two existing multigrid methods, and then investigate more
details on the proposed preconditioner $B$.
\vskip 0.1in
\noindent$\bullet$ Comparisons with the standard multigrid method with overlapping Schwarz smoothers

The preconditioner $B$ defined in the previous two subsection looks like the standard multigrid preconditioner with overlapping Schwarz smoothers,
but the two preconditioners have essential differences. In order to explain the differences in details, we first describe this standard preconditioner for the current situation.

As in Subsection 3.1, let $\Omega$ be decomposed into the union of several quasi-uniform and regular coarse elements with the size $h_0$,
where each coarse element is just the union of some fine elements in ${\mathcal T}_h$.  Let
${\mathcal T}_{h_0}$ denote the resulting partition, i.e., the set of all the coarse elements. For every element in ${\mathcal T}_{h_0}$, we continue
such decomposition and obtain several quasi-uniform and regular coarse elements with the size $h_1<h_0$. The resulting partition is denoted by $\hat{\mathcal T}_{h_1}$.
As usual, we repeat the above decomposition process and we can build refining finite element partitions: ${\mathcal T}_{h_0}$, $\hat{\mathcal T}_{h_1},\cdots,
\hat{\mathcal T}_{h_J}$ with the sizes $h_0$, $h_1,\cdots, h_J$ satisfying $h<h_J<\cdots <h_1<h_0$. For $j=1,\cdots,J$, let $V_p(\hat{\mathcal T}_{h_j})$ denote
the plane wave finite element space associated with the finite element partitions $\hat{\mathcal T}_{h_j}$. Then we obtain
the multilevel space decomposition
\begin{equation}
V_p({\mathcal T}_{h})=V_p({\mathcal T}_{h_0})+\sum\limits_{j=1}^JV_p(\hat{\mathcal T}_{h_j}).\label{6.decomposition1}
\end{equation}

In order to define overlapping Schwarz smoothers, we decompose the space $V_p(\hat{\mathcal T}_{h_j})$ ($j\geq 1$) into the sum of smaller subspaces.
For each $D\in \hat{\mathcal T}_{h_j}$, we enlarge $D$ with the thickness of one (coarse) element to a larger domain $\tilde{D}$, i.e., $\tilde{D}$ is the union of $D$
and the coarse elements adjoining $D$, where the added elements belong to $\hat{\mathcal T}_{h_j}$. Then $D$ and the added elements constitute a coarse finite element partition of $\tilde{D}$
, which is denoted by $\hat{\mathcal T}^{\tilde{D}}_{h_j}$. Let $V_p(\hat{\mathcal T}^{\tilde{D}}_{h_j})$ denote the plane wave finite element space associated with the finite element partition
$\hat{\mathcal T}^{\tilde{D}}_{h_j}$, i.e., the restriction of $V_p(\hat{\mathcal T}_{h_j})$ on the subdomain $\tilde{D}$. Then we have the ``overlapping" space decomposition of the $j$th-level
coarse space
\begin{equation}
V_p(\hat{\mathcal T}_{h_j})=\sum_{D\in \hat{\mathcal T}_{h_j}}V_p(\hat{\mathcal T}^{\tilde{D}}_{h_j})\quad\quad(j\geq 1). \label{6.decomposition0}
\end{equation}
Combing this decomposition with (\ref{6.decomposition1}), gives the new multilevel decomposition of the original space
\begin{equation}
V_p({\mathcal T}_{h})=V_p({\mathcal T}_{h_0})+\sum\limits_{j=1}^J\sum_{D\in \hat{\mathcal T}_{h_j}}V_p(\hat{\mathcal T}^{\tilde{D}}_{h_j}).\label{6.decomposition2}
\end{equation}

As in Subsection 3.2, let $A_0$ be the coarsest solver associated with $V_p({\mathcal T}_{h_0})$. We define $\hat{A}^{\tilde{D}}_{h_j}: V_p(\hat{\mathcal T}^{\tilde{D}}_{h_j}):\rightarrow V_p(\hat{\mathcal T}^{\tilde{D}}_{h_j})$
as the restriction of $A$ on $V_p(\hat{\mathcal T}^{\tilde{D}}_{h_j})$, and use $\hat{Q}^{\tilde{D}}_{h_j}: V_p(\hat{\mathcal T}_{h_j}):\rightarrow V_p(\hat{\mathcal T}^{\tilde{D}}_{h_j})$ to denote the $L^2$ projector.
Then we define $j$th-level solver
$$ \hat{B}^{-1}_j=\sum_{D\in \hat{\mathcal T}_{h_j}}(\hat{A}^{\tilde{D}}_{h_j})^{-1}\hat{Q}^{\tilde{D}}_{h_j}\quad\quad(1\leq j\leq J) $$
and the preconditioner
\begin{equation}
\hat{B}^{-1}=A_0^{-1}Q_0+\sum\limits_{j=1}^J\hat{B}^{-1}_jQ_j.\label{6.preconditioner1}
\end{equation}
For $j\geq 1$, the operator $\hat{B}^{-1}_j$ is called the overlapping Schwarz smoother at $j$th-level, and the preconditioner $\hat{B}$ is called the
multigrid preconditioner with overlapping Schwarz smoothers (MG-Schwarz). For each coarse element $D\in\hat{\mathcal T}_{h_j}$, if we do not enlarge $D$ into the larger subdomain $\tilde{D}$ and replace
the subspace $V_p(\hat{\mathcal T}^{\tilde{D}}_{h_j})$ in (\ref{6.decomposition2}) by $V_p(\hat{\mathcal T}^{D}_{h_j})$ itself, then the corresponding preconditioner $\hat{B}$ is just the multigrid preconditioner
with Jacobi smoothers (MG-Jacobi). In applications, the action of the smoother $\hat{B}^{-1}_j$ may be repeated several times by
Richardson iterations. Notice that we have not considered the more general situation, in which the subdomain $\tilde{D}$ contains more elements for each $D\in\hat{\mathcal T}_{h_j}$,
since the implementation of the resulting smoothers has greater cost.

Now we give some comparisons between the preconditioner $B$ defined in (\ref{3.preconditioner}) and the preconditioner $\hat{B}$ defined in (\ref{6.preconditioner1}). We need only to
compare the two multilevel space decompositions (\ref{decom-new}) and (\ref{6.decomposition2}).

\noindent {\bf Similarity}: for both multilevel space decompositions, the subspaces in each level (except the coarsest level) are overlapping each other.

\noindent {\bf Differences}:

(1) the two space decompositions are constructed in different ways. For the space decomposition (\ref{6.decomposition2}), we first have the multilevel decomposition (\ref{6.decomposition1}), and then construct independently the overlapping decomposition (\ref{6.decomposition0}) for each level coarse space. However, for the space decompositions (\ref{decom-new}), we first construct the overlapping decomposition (\ref{3.local_decomposition}) on each ``enlarged subdomain", and then use all these local overlapping decompositions to derive recursively the global multilevel space decomposition (\ref{decom-new}).

(2) the two space decompositions have different structures. The design of the overlapping decomposition (\ref{6.decomposition0}) only changes the structure of $j$th-level space itself, but does not improve the relation of the coarse spaces at different levels. This means that the structure of the space decomposition (\ref{6.decomposition2}) has no essential difference from that in the multilevel preconditioner with Jacobi smoothers. From the construction of the space decomposition (\ref{decom-new}), we know that the space decomposition (\ref{decom-new}) locally possesses the structure of the space decomposition in the overlapping domain decomposition method, and so
the overlapping subspaces $V_p({\mathcal T}^{D}_{d_j})$ at different levels have inherent connections. It is easy to see that the space decomposition (\ref{decom-new}) is independent of the space decomposition in the standard multigrid preconditioners. Some comparison results for them will be given in Table \ref{table-new3} of Section 6.
\vskip 0.1in
\noindent$\bullet$ Comparisons with the wave-ray multigrid methods

The wave-ray multigrid methods (see \cite{LivBrand2006} and \cite{LeeMccor2000}) were designed for solving Helmholtz system generated by the discretization with finite difference or the nodal finite elements.

As in the first part of this subsection, let $h_j$ denote the size of the coarse elements at $j$-th level. It is well known that, when $h_j$ is relatively large comparing the value of $1/\omega$, the oscillatory
error components at $j$-th level can not be efficiently reduced by the standard multigrid methods. The basic idea of the wave-ray multigrid methods is to approximate such oscillatory error components at $j$-th level
by the following functions
$$ w^j({\x})=\sum_{l=1}^{L_j}a^j_l({\bf x})e^{i\omega(\boldsymbol{\alpha_l}\cdot{\bf x})}, $$
where $a^j_l({\bf x})$ are smooth functions, which are called {\it ray envelope functions} in \cite{LivBrand2006}; the wave direction vectors
$\{\boldsymbol{\alpha_l}\}$ may be different from that given in Subsection 2.4. The number $L_j$ of the wave directions increases when the value $\omega h_j$ increases.

Since the original error components can not be directly expressed as the form of the function $w^j({\bf x})$, some exponential interpolations need to be constructed (see \cite{LeeMccor2000}).
These exponential interpolations were defined by the Fourier components (ray elements, plane wave functions) $e^{i\omega(\boldsymbol{\alpha_l}\cdot{\bf x})}$, and were used to achieve a transformation between
the original error components and the ray envelope functions. In the wave-ray multigrid methods, the approximation of oscillatory error components was transformed into the approximation of smooth ray envelope
functions by using the exponential interpolations. While the smooth ray envelope functions can be approximated by the standard multigrid methods. Then the oscillatory error components can be reduced on relatively
coarse girds. The implementation of the wave-ray multigrid methods involves many technical details, for example, how to choose suitable wave direction vectors $\{\boldsymbol{\alpha_l}\}$. The cost in the wave-ray
multigrid methods depends on the value of the wave number $L_j$ and the calculation of the exponential interpolations.

Notice that both the wave-ray multigrid method and the multilevel method introduced in this paper are based on the plane wave functions $e^{i\omega(\boldsymbol{\alpha_l}\cdot{\bf x})}$, in essence, use
the ``good" approximate property of the plane wave functions for oscillatory solutions. However, the roles of the plane wave functions are different in the two kinds of methods: the plane wave functions are used to
define discretization basis functions in this paper; while, the plane wave functions are
only auxiliary {\it weight} functions in the wave-ray multigrid methods. As to the multilevel methods themselves, the multilevel method described in the previous two subsections has no
relation with the wave-ray multigrid methods, since the wave-ray multigrid methods still use the standard multigrid framework to approximate the smooth ray envelope functions.

\vskip 0.1in

\noindent$\bullet$ On the efficiency of the proposed multilevel method.

In Section 6, we will test several examples to illustrate the efficiency of the proposed multilevel preconditioners (including some variants of $B$, see Sections 4-5).
Besides, we will give numerical comparisons among the proposed multilevel preconditioner $B$, the multigrid preconditioner $\hat{B}$ with overlapping Schwarz smoothers and the multigrid preconditioner
with Jacobi smoothers. As we will see, the multilevel preconditioner $B$ designed in the previous two subsections is robust even for large $\omega$. However,
the multilevel preconditioner $\hat{B}$ with overlapping Schwarz smoothers can only slightly improve the convergence rate of the multilevel preconditioner with
Jacobi smoothers. In this part,we try to give some explanations
to the effectiveness of the preconditioner $B$.

The first reason is that the plane wave functions can approximate the oscillatory solution of the Helmholtz equation very well (which is just the motive of the wave-ray multigrid methods), but it is not the unique reason
of the effectiveness. In fact, if we decrease the thickness of the overlap to be one fine element in the overlapping space decomposition (\ref{decom-new}), then the resulting multilevel preconditioner
has almost the same convergence rate with the multigrid preconditioner $\hat{B}$ with Schwarz smoothers (see the results reported in Table \ref{table-new2} of Section 6).
The second reason is that the space decomposition (\ref{decom-new}) possesses ``good" structure, as explained in the first part of this subsection. We would like to explain this point more clearly.
It is known that an overlapping domain decomposition preconditioner with several subdomains only is always stable even for the Helmholtz equations with large wave numbers (see the results listed in Table \ref{table-new1} of
Section 6). Thus, since the number $N_0$ of overlapping subdomains is fixed and not large, the overlapping decomposition (\ref{3.local_decomposition}) (and (\ref{3.intial_decomposition})) is stable for each $D\in{\mathcal S}_{j-1}$ even for large $\omega$. This means that the global space decomposition (\ref{decom-new}), which is defined by the local space decompositions (\ref{3.intial_decomposition}) and (\ref{3.local_decomposition}), should be also stable even for large $\omega$. Notice that each local space $V_p({\mathcal T}^{\tilde{D}_r}_h)$ has too high dimension unless $j$ is large, so we have to make multilevel
decomposition.

\vskip 0.1in

\noindent$\bullet$ Computational cost for the implementation of the proposed preconditioner $B$.

In applications, the action of $B^{-1}$ is implemented {\it in parallel}. Thus we should not investigate the computational complexity
for the implementation of $B^{-1}$ as successive algorithm. But, for completeness, we still
estimate the computational complexity in the usual way.

As in Section 2, let $N$ denote the number of the fine elements in ${\mathcal T}_h$. It is easy to see that the numbers of different subproblems needed to be solved in {\bf Algorithm 3.2} and {\bf Algorithm 3.3}
are not greater than~~$min\{N_0^{j-1},~~N\}$ $(j=1,\cdots,J)$ and~~$min\{N_0^J,~~N\}$, respectively.
If we require that the number of the fine elements contained in each
$K\in {\mathcal S}_J$ almost equals $N_0$ (refer to Remark 3.1), we can verify that the level number $J$ should be $c_0(\log_{N_0} N)$, where $c_0$ is a positive constant depending on $N_0$ and $\theta_0$.
Notice that each subproblem to be solved has $N_0p$ unknowns, so its solution has the computational cost $O((N_0p)^3)$. Then the computational complexity for the implementation of $B^{-1}$
can be estimated as follows
$$ {\mathcal N}_{cost}\leq C(N_0p)^3J \min\{N_0^J,N\}\leq CN_0^3p^3(\log_{N_0} N)N=CN_0^3p^2(\log_{N_0} N)(pN). $$
Then we have
$$ {\mathcal N}_{cost}\leq C N_0^3p^2(\log_{N_0} N) N_{dof}, $$
where $N_{dof}=pN$ denotes the dimension of the original fine grid system (\ref{eq34}). Since $N_0$ is a constant, the
computational cost is estimated by ${\mathcal N}_{cost}=O(p^2(\log_{N_0} N) N_{dof}).$  This means that, even if we implement the action of $B^{-1}$ in {\it successive}
manner, the resulting computational complexity is almost the optimal. Since the solution of each subproblem has very small cost $O((pN_0)^3)$,
the preconditioner $B$ implemented in parallel should be much cheaper than the direct solvers.

\section{A multilevel overlapping preconditioner with smoothers}
In this section, we design an improvement of the preconditioner $B$ to further reduce the cost for implementing
the solvers $B^{-1}_j$ ($j=1,\cdots,J$) and $\tilde{B}^{-1}_J$ described in {\bf Algorithm 3.2} and {\bf Algorithm 3.3}. The basic idea
is to replace the solvers $B^{-1}_j$ ($j=1,\cdots,J$) and $\tilde{B}^{-1}_J$ by Jacobi-type smoothers. To this end, we first give exact definitions of the smoothers.

For $j=1,\cdots,J$ and $D\in {\mathcal S}_{j-1}$, let $V_p({\mathcal T}^{D}_{d_j})$ denote the local coarse space defined in Subsection 3.1. We want to further decompose
each space $V_p({\mathcal T}^{D}_{d_j})$ into the sum of several smaller spaces. Notice that the support set of the functions in $V_p({\mathcal T}^{D}_{d_j})$
is $D$, which is the union of $N_0$ coarse elements $D_1,\cdots,D_{N_0}$ in ${\mathcal T}^{D}_{d_j}$. Thus we need only to define subspaces on the coarse elements.

As in Subsection 3.1, let ${\mathcal Q}_p$ denote the space of $p$ plane wave shape functions. For a coarse element $D_r$ in ${\mathcal T}^{D}_{d_j}$, define
$$
V_p(D_r)=span\{v\in L^2(\Omega):~v|_{D_r}\in {\mathcal Q}_p;~~supp~v\subset D_r
\}=\{v\in V_p({\mathcal T}^{D}_{d_j}):~supp~v\subset D_r\}. $$
$$ (j=1,\cdots,J; D\in {\mathcal S}_{j-1};~r=1,\cdots,N_0)$$
In other words, $V_p(D_r)$ is the restriction space of $V_p({\mathcal
T}^{D}_{d_j})$ on $D_r\subset D$. It is clear that
the space $V_p(D_r)$ has the dimension $p$ (but the dimension of $V_p({\mathcal T}^{D}_{d_j})$ equals $N_0p$). Then
$$ V_p({\mathcal T}^{D}_{d_j})=\sum_{r=1}^{N_0} V_p(D_r), $$
and so
$$ V_p({\mathcal T}_{d_j})=\sum_{D\in {\mathcal S}_{j-1}}\sum_{r=1}^{N_0}V_p(D_r). $$

Similarly, for each {\it fine} element $E\in \tilde{\mathcal T}^K_h$, define
$$ V_p(E)=\{v\in V_p(\tilde{\mathcal T}^K_h):~~supp~ v\subset E\}\quad\quad(K\in {\mathcal S}_{J};~E\in \tilde{\mathcal T}^K_h). $$
It is clear that the dimension of $V_p(E)$ equals $p$ and we have
$$ V_p(\tilde{\mathcal T}^K_h)=\sum_{E\in\tilde{\mathcal T}^K_h}V_p(E). $$
Then
$$ V_p(\tilde{\mathcal T}^J_h)=\sum_{K\in {\mathcal S}_{J}}\sum_{E\in\tilde{\mathcal T}^K_h}V_p(E). $$

Based on the above space decompositions, we can define Jacobi-type smoothers in the natural manner.

Let $m_0$ be a given positive integer.
 The desired smoothers $R^{(m_0)}_j$ ($j=1,\cdots,J$) and $\tilde{R}^{(m_0)}_J$ are defined by the following algorithms.

{\bf Algorithm 4.1}. For $\eta\in V_p({\mathcal T}_{d_j})$, the function $w_{\eta}=(R^{(m_0)}_j)^{-1}\eta\in V_p({\mathcal T}_{d_j})$ can be obtained as follows:

Step 1. Let $w^0\in V_p({\mathcal T}_{d_j})$ be an initial guess. Assume that $w^{l-1}$ ($l=1,\cdots,m_0$) has been gotten. For $D\in{\mathcal S}_{j-1}$ and elements $D_r\in {\mathcal T}^{D}_{d_j}$, computing $\hat{w}^l_{D_r}\in V_p(D_r)$ in parallel by
\[
a(\hat{w}^l_{D_r},v)=(\eta,v)-a(w^{l-1},v),\quad \forall v\in V_p(D_r),
\]
and set
$$ w^l=w^{l-1}+\sum_{D\in{\mathcal S}_{j-1}}\sum_{r=1}^{N_0}\hat{w}^l_{D_r}\quad (l=1,\cdots,m_0); $$

Step 2. Define $w_{\eta}=w^{m_0}$.

\vskip 0.2in

{\bf Algorithm 4.2}. For $\eta\in V_p(\tilde{\mathcal T}^J_h)$, the function $w_{\eta}=(\tilde{R}^{(m_0)}_J)^{-1}\eta\in V_p(\tilde{\mathcal T}^J_h)$ can be obtained as follows:

Step 1. Let $w^0\in V_p(\tilde{\mathcal T}^J_h)$ be an initial guess. Assume that $w^{l-1}$ ($l=1,\cdots,m_0$) has been gotten. For $K\in {\mathcal S}_J$ and {\it fine} elements $E\in\tilde{\mathcal T}^K_h$, computing $\hat{w}^l_E\in V_p(E)$ in parallel by
\[
a(\hat{w}^l_E, v)=(\eta,v)-a(w^{l-1},v),\quad \forall v\in V_p(E),
\]
and set
$$ w^l=w^{l-1}+\sum_{K\in{\mathcal S}_{J}}\sum_{E\in\tilde{\mathcal T}^K_h}\hat{w}^l_E\quad (l=1,\cdots,m_0); $$

Step 2. Define $w_{\eta}=w^{m_0}$.

\vskip 0.2in

Next we define a new multilevel preconditioner.

Let $m_0$ be a given positive integer, and let $(R^{(m_0)}_j)^{-1}$ and $(\tilde{R}^{(m_0)}_J)^{-1}$ denote the smoothers defined by {\bf Algorithm 4.1} and {\bf Algorithm 4.2}, respectively. Define the additive preconditioner
$$ (B_s^{(m_0)})^{-1}=A^{-1}_0Q_0+\sum_{j=1}^{J}(R^{(m_0)}_j)^{-1}Q_j+(\tilde{R}^{(m_0)}_J)^{-1}\tilde{Q}_J.$$

In applications, we can choose the positive integer $m_0$ as $m_0=2,3$. The action of $(B_s^{(m_0)})^{-1}$ can be implemented as in {\bf Algorithm 3.1}, provided that the solvers $B_j$ ($j=1,\cdots,J$) and $\tilde{B}_J$ are replaced with $R^{(m_0)}_j$ ($j=1,\cdots,J$) and $\tilde{R}^{(m_0)}_J$ defined by {\bf Algorithm 4.1} and {\bf Algorithm 4.2}. Since the actions of $(R^{(m_0)}_j)^{-1}$ ($j=1,\cdots,J$) and $(\tilde{R}^{(m_0)}_J)^{-1}$ are implemented in smaller spaces, one of which is defined on an (coarse or fine) element and has only $p$ degree of freedoms, the preconditioner $B_s^{(m_0)}$ is cheaper than the preconditioner $B$.
Numerical experiments in Section 6 will indicate that the new variant has faster convergence than the preconditioner $B$.

\begin{remark}
According to the discussions in Subsection 3.3, the preconditioner $B_s^{(m_0)}$ is different from the standard multigrid preconditioners since the space decomposition defining $B_s^{(m_0)}$ possesses
different structure from the one corresponding to the standard multigrid preconditioners. The differences between $B_s^{(m_0)}$ and the multigrid preconditioner $\hat{B}$ with overlapping Schwarz smoothers
are more obvious: each subproblem (except the coarsest problem) to be solved in $B_s^{(m_0)}$ has $p$ unknowns only, but each subproblem to be solved in $\hat{B}$ has $n_{\tilde{D}}\times p$ unknowns, where
$n_{\tilde{D}}$ denotes the number of the (coarse) elements contained in a subdomain ${\tilde{D}}$. The proposed method is not called as multigrid method, since the sets ${\mathcal T}_{d_j}$ and $\tilde{\mathcal T}_h^J$
defining the multilevel spaces do not constitute grids on $\Omega$ yet. For convenience, we called the preconditioner $B_s^{(m_0)}$ as {\it multilevel overlapping preconditioners with smoothers} (MOPS).
\end{remark}

\begin{remark} Notice that the dimension of the coarsest space $V_p({\mathcal T}^{\Omega}_{d_0})$ equals $N_0p$, with
$N_0$ being a constant independent of $\omega$, $h$ and $p$.
Thus, it is cheap to realize the action of $A_0^{-1}$ appearing in the preconditioner $B_s^{(m_0)}$ by the direct method (the values of $N_0$ and $p$ are not large). Of course, the action of $A_0^{-1}$ can be also replaced
by implementing a cheaper preconditioner of $A_0$. It is easy to construct such a cheaper preconditioner for $A_0$ since the space $V_p({\mathcal T}^{\Omega}_{d_0})$ is defined on $N_0$ coarsest elements with fixed size $d_0$.
\end{remark}

\section{Multiplicative variants of the preconditioner $B_s^{(m_0)}$}
In this section, we design several multiplicative multilevel preconditioners to accelerate the convergence of the additive preconditioner $B_s^{(m_0)}$.
\subsection{A basic multiplicative preconditioner}

In this subsection, we introduce a simple multiplicative preconditioner.

Define the operator
\[
P_0=A_0^{-1}Q_0A.
\]
Then $P_0$ is the energy projector from $V_p(\mathcal{T}_h)$ into the coarsest space $V_p(\mathcal{T}_{d_0})$.
Let $R^{(m_0)}_j$ ($j=1,\cdots,J$) and $\tilde{R}^{(m_0)}_J$ be the smoothers defined in the last section,
and set
$$ T^{(m_0)}_j=(R^{(m_0)}_j)^{-1}Q_jA\quad(j=1,\cdots,J)\quad\mbox{and}\quad \tilde{T}^{(m_0)}_J=(\tilde{R}^{(m_0)}_J)^{-1}\tilde{Q}_JA. $$
Let $I$ denote the identity operator on $V_p(\mathcal{T}_h)$. Associated with the space decomposition
(\ref{decom6}), a multiplicative variant of $B_s^{(m_0)}$ is defined by
$$ (M^{(m_0)}_1)^{-1}=\big(I-(I-P_0)(I-T^{(m_0)}_1)\cdots(I-T^{(m_0)}_{J})(I-\tilde{T}^{(m_0)}_J)\big)A^{-1}.$$
The error propagation operator of $M^{(m_0)}_1$ is
$$ I-(M_1^{(m_0)})^{-1}A=(I-P_0)(I-T^{(m_0)}_1)\cdots(I-T^{(m_0)}_{J})(I-\tilde{T}^{(m_0)}_J). $$

The action of $(M^{(m_0)}_1)^{-1}$ can be described by the following algorithm.

\vskip 0.2in

{\bf Algorithm 5.1}. For $\xi\in V_p(\mathcal{T}_h)$, the function $u_\xi=(M^{(m_0)}_1)^{-1}\xi\in V_p(\mathcal{T}_h)$ can be obtained as follows:

Step 1. Computing $\tilde{u}^J_h\in V_p(\tilde{\mathcal T}^{J}_h)$ by
\[
(\tilde{R}^{(m_0)}_J\tilde{u}^J_h,v_h)=(\xi,v_h),\quad \forall v_h\in V_p(\tilde{\mathcal T}^{J}_h);
\]

Step 2. Computing $u_{d_{J}}\in V_p({\mathcal T}_{d_{J}})$ by
\[
(R^{(m_0)}_{J}u_{d_{J}},v)=(\xi,v)-a(\tilde{u}^J_h,v),\quad \forall v\in V_p({\mathcal T}_{d_{J}}),
\]
and set $u_{J}=\tilde{u}^J_h+u_{d_{J}}$;

Step 3. Let $j=J,\cdots,2$. If we have obtained $u_j\in V_p(\mathcal{T}_h)$, then compute $u_{d_{j-1}}\in V_p(\mathcal{T}_{d_{j-1}})$ by
\[
(R^{(m_0)}_{j-1}u_{d_{j-1}},v)=(\xi,v)-a(u_j,v),\quad \forall v\in V_p(\mathcal{T}_{d_{j-1}}),
\]
and set
\[
u_{j-1}=u_{j}+u_{d_{j-1}}\quad (j=J,\cdots,2);
\]

Step 4. Computing $u_0\in V_p(\mathcal{T}_{d_0})$ by
\[
(A_0u_0,v)=(\xi,v)-a(u_1,v),\quad \forall v\in V_p(\mathcal{T}_{d_0});
\]

Step 5. Set
\[
u_\xi=u_1+u_0.
\]

\subsection{The standardly symmetrized multiplicative preconditioner}
In this subsection we consider the case of PWLS method. Then the operator $A$ is Hermitian positive definite
with respect to the inner product $(A\cdot,\cdot)$. Thus, we need to define a symmetrization of the
preconditioner $M^{(m_0)}_1$.

The standardly symmetrized preconditioner of $M^{(m_0)}_1$ is defined as
$$
(M^{(m_0)}_2)^{-1}=\big(I-(I-\tilde{T}^{(m_0)}_J)(I-T^{(m_0)}_{J})\cdots(I-T^{(m_0)}_1)(I-P_{0})(I-T^{(m_0)}_1)
\cdots(I-T^{(m_0)}_{J})(I-\tilde{T}^{(m_0)}_J)\big)A^{-1}. $$
The error propagation operator of $M^{(m_0)}_2$ is
$$ I-(M^{(m_0)}_2)^{-1}A=(I-\tilde{T}^{(m_0)}_J)(I-T^{(m_0)}_{J})\cdots(I-T^{(m_0)}_1)(I-P_{0})(I-T^{(m_0)}_1)
\cdots(I-T^{(m_0)}_{J})(I-\tilde{T}^{(m_0)}_J).$$
For the case of PWLS method, the operators $T^{(m_0)}_j$ and $\tilde{T}^{(m_0)}_J$ are Hermitian positive definite with respect to the inner product $(A\cdot,\cdot)$. As a result, the operator $(M^{(m_0)}_2)^{-1}$ is also Hermitian and positive definite with respect to the same inner product.

The action of $(M^{(m_0)}_2)^{-1}$ can be described by the following algorithm.

\vskip 0.2in

{\bf Algorithm 5.2}. For $\xi\in V_p(\mathcal{T}_h)$, the function $u_\xi=(M^{(m_0)}_2)^{-1}\xi\in V_p(\mathcal{T}_h)$ can be obtained as follows:

Step 1. Computing $\tilde{w}^J_h\in V_p(\tilde{\mathcal T}^{J}_h)$ by
\[
(\tilde{R}^{(m_0)}_J\tilde{w}^J_h,v_h)=(\xi,v_h),\quad \forall v_h\in V_p(\tilde{\mathcal T}^{J}_h);
\]

Step 2. Computing $w_{d_{J}}\in V_p({\mathcal T}_{d_{J}})$ by
\[
(R^{(m_0)}_{J}w_{d_{J}},v)=(\xi,v)-a(\tilde{w}^J_h,v),\quad \forall v\in V_p({\mathcal T}_{d_{J}}),
\]
and set $w_{J}=\tilde{w}^J_h+w_{d_{J}}$;

Step 3. Let $j=J,\cdots,2$. If we have obtained $w_j\in V_p(\mathcal{T}_h)$, then compute $w_{d_{j-1}}\in V_p(\mathcal{T}_{d_{j-1}})$ by
\[
(R^{(m_0)}_{j-1}w_{d_{j-1}},v)=(\xi,v)-a(w_j,v),\quad \forall v\in V_p(\mathcal{T}_{d_{j-1}}),
\]
and set
\[
w_{j-1}=w_j+w_{d_{j-1}}\quad (j=J,\cdots,2);
\]

Step 4. Computing $u_{d_0}\in V_p(\mathcal{T}_{d_0})$ by
\[
(A_0u_{d_0},v)=(\xi,v)-a(w_1,v),\quad \forall v\in V_p(\mathcal{T}_{d_0}),
\]
and set $u_0=w_1+u_{d_0}$;

Step 5. Let $j=1,\cdots,J$. If we have obtained $u_{j-1}\in V_p(\mathcal{T}_h)$, then compute $u_{d_j}\in V_p(\mathcal{T}_{d_j})$ by
\[
(R^{(m_0)}_ju_{d_j},v)=(\xi,v)-a(u_{j-1},v),\quad \forall v\in V_p(\mathcal{T}_{d_j}),
\]
and set
\[
u_j=u_{j-1}+u_{d_j}\quad (j=1,\cdots,J).
\]

Step 6. Computing $\tilde{u}^J_h\in V_p(\tilde{\mathcal T}^{J}_h)$ by
\[
(\tilde{R}^{(m_0)}_J\tilde{u}^J_h,v_h)=(g,v_h)-a(u_{J},v_h),\quad \forall v_h\in V_p(\tilde{\mathcal T}^{J}_h);
\]

Step 7. Set
\[
u_\xi=u_{J}+\tilde{u}^J_h.
\]

\subsection{A non-standard symmetrized multiplicative preconditioner}

In this subsection, we still consider the case of PWLS method. Define the operator $T^{(m_0)}: V_p(\mathcal{T}_h)\to V_p(\mathcal{T}_h)$
by
$$ T^{(m_0)}=I-(I-T^{(m_0)}_1)\cdots(I-T^{(m_0)}_{J})(I-\tilde{T}^{(m_0)}_J)
(I-T^{(m_0)}_{J})\cdots(I-T^{(m_0)}_1). $$
Then $T^{(m_0)}$ is Hermitian positive definite with respect to the inner product $(A\cdot,\cdot)$. A non-standard symmetrized preconditioner of $M^{(m_0)}_1$ can be defined as (refer to \cite{HuL2002})
$$ (M^{(m_0)}_3)^{-1}=\big(I-(I-P_0)(I-T^{(m_0)})\big)A^{-1}\quad\quad(\mbox{the}~~PWLS~~\mbox{method}). $$

It can be verified that the restriction of $M^{(m_0)}_3$ on $(V_p({\mathcal T}_{d_0}))^{\bot}$ is Hermitian positive definite with respect to the inner product $(A\cdot,\cdot)$ (refer to \cite{HuL2002}). The error propagation operator of $M^{(m_0)}_3$ is
$$ I-(M^{(m_0)})^{-1}_3A=(I-P_0)(I-T). $$

The action of $(M^{(m_0)}_3)^{-1}$ can be described by the following algorithm.

\vskip 0.2in

{\bf Algorithm 5.3}. For $\xi\in V_p(\mathcal{T}_h)$, the function $u_\xi=(M^{(m_0)}_3)^{-1}\xi\in V_p(\mathcal{T}_h)$ can be obtained as follows:

Step 1. Computing $w_1\in V_p(\mathcal{T}_{d_1})$ by
\[
(R^{(m_0)}_1w_1,v)=(\xi,v),\quad \forall v\in V_p(\mathcal{T}_{d_1});
\]

Step 2. Let $j=2,\cdots,J$. If we have obtained $w_{j-1}\in V_p(\mathcal{T}_h)$, then compute $w_{d_j}\in V_p(\mathcal{T}_{d_j})$ by
\[
(R^{(m_0)}_jw_{d_j},v)=(\xi,v)-a(w_{j-1},v),\quad \forall v\in V_p(\mathcal{T}_{d_j}),
\]
and set
\[
w_j=w_{j-1}+w_{d_j}\quad (j=2,\cdots,J);
\]

Step 3. Computing $\tilde{w}^J_h\in V_p(\tilde{\mathcal T}^{J}_h)$ by
\[
(\tilde{R}^{(m_0)}_J\tilde{w}^J_h,v_h)=(\xi,v_h)-a(w_{J},v_h),\quad \forall v_h\in V_p(\tilde{\mathcal T}^{J}_h),
\]
and set $\tilde{w}_J=w_{J}+\tilde{w}^J_h$;

Step 4. Computing $u_{d_{J}}\in V_p({\mathcal T}_{d_{J}})$ by
\[
(R^{(m_0)}_{J}u_{d_{J}},v)=(\xi,v)-a(\tilde{w}_J,v),\quad \forall v\in V_p({\mathcal T}_{d_{J}}),
\]
and set $u_{J}=\tilde{w}_J+u_{d_{J}}$;

Step 5. Let $j=J,\cdots,2$. If we have obtained $u_j\in V_p(\mathcal{T}_h)$, then compute $u_{d_{j-1}}\in V_p(\mathcal{T}_{d_{j-1}})$ by
\[
(R^{(m_0)}_{j-1}u_{d_{j-1}},v)=(\xi,v)-a(u_j,v),\quad \forall v\in V_p(\mathcal{T}_{d_{j-1}}),
\]
and set
\[
u_{j-1}=u_{j}+u_{d_{j-1}}\quad (j=J,\cdots,2);
\]

Step 6. Computing $u_0\in V_p(\mathcal{T}_{d_0})$ by
\[
(A_0u_0,v)=(\xi,v)-a(u_1,v),\quad \forall v\in V_p(\mathcal{T}_{d_0});
\]

Step 7. Set
\[
u_\xi=u_1+u_0.
\]

\begin{remark} The actions of $(R^{(m_0)}_j)^{-1}$ ($j=1,\cdots,J$) and $(\tilde{R}^{(m_0)}_J)$ used in {\bf Algorithm 5.1}
-{\bf Algorithm 5.3} are implemented by {\bf Algorithm 4.1} and {\bf Algorithm 4.2}, respectively. Notice that the solver $(\tilde{R}^{(m_0)}_J)^{-1}$ is implemented only one time in {\bf Algorithm 5.3} (such solver needs to be implemented for two times in {\bf Algorithm 5.2}), so the preconditioner $M^{(m_0)}_3$ is cheaper than $M^{(m_0)}_2$. It is interesting that the numerical results reported in Section 6 indicate that $M^{(m_0)}_3$ has faster convergence than $M^{(m_0)}_2$ (some explanations to the kind of phenomenon have been given in \cite{HuL2002}).
\end{remark}

\section{Numerical experiments}
In this section we report numerical results to illustrate that the new preconditioners are effective for solving Helmholtz equations with large wave numbers.

In the examples tested in this part, we choose $\Omega$ as the rectangle $[0,2]\times[0,1]$, and we adopt a uniform
partition $\mathcal {T}_h$ for the domain $\Omega$ as follows:
$\Omega$ is divided into some small rectangles with the
same size, where $h$ denotes the length of the longest edge of the
elements. Let $n_h$ denote the number of elements generated by the partition $\mathcal {T}_h$, and let $p$ denote the number of plane wave basis functions in one element.
Then the dimension of the original fine grid system (\ref{eq34}) is $N_{dof}=n_h\times p$.

We choose the mesh size $h$ and the number $p$ of plane wave basis functions in one element according to the following rule:  when the wave numbers increase, the scale of the discrete problem is increased (either $h$ decreases or $p$ increases) in a suitable manner such that accepted relative $L^2$ errors of the approximation can be kept.
In the numerical experiments below, we choose $h\approx 2/\omega$ and slightly increase $p$ when $\omega$ increases.

We need to give a rule for the multilevel overlapping domain decomposition. For convenience, we consider only an easily implemented rule, i.e., the overlap degree $\theta_0=1$, for the main experiments.
Let $\Omega$ be divided into $2^n\times 2^n$ ($n\geq 3$) rectangle elements with the same size. We divide $\Omega$ into 4 parts in each direction ($x$-coordinate axis direction or $y$-coordinate axis direction) to build the coarsest partition ${\mathcal T}_{d_0}$, with $d_0$ being a constant independent of the wave number $\omega$ and the fine mesh size $h$. This means that the coarsest partition contains $4\times4$ (coarse) rectangular elements with the same size, and so $N_0=16$. Define the enlarged subdomain of each (coarse) element as the union of the (coarse) element itself and its neighboring (coarse) elements, where the definition of the enlarged subdomain was given in Subsection 3.1. We repeat the above process to decompose each enlarged subdomain into $4\times 4$ rectangles, but the rectangles may have different sizes since the number of the elements contained in a enlarged rectangle may be not divisible by $16$. For this case, we still divide the enlarged rectangle into 4 parts in each direction such that the number of elements in each part is almost the same. We continue the above process, and the decomposition stops when the number of elements in each enlarged subdomain associated with the current level is less than $5\times 5$.

Throughout this section, we always use $B$, $B_s^{(m_0)}$ and $M^{(m_0)}_l$ ($l=1,2,3$) to denote the proposed multilevel preconditioners with the above decomposition rule.

For the PWLS method, we set $\alpha=\omega^2$ and $\beta=1$; for the PWDG method, we set $\alpha=\beta=\delta={1\over 2}$.
Since the stiffness matrix of PWLS method is Hermitian positive definite, we can solve the system by PCG method. While the stiffness matrix of PWDG is not Hermitian, we solve it by PGMRES method. For one iterative step, PCG method is cheaper than PGMRES method. The stopping criterion in the iterative algorithms is that the
relative $L^2$-norm $\epsilon$ of the residual of the iterative
approximation satisfies $\epsilon < 1.0e-6 \ $.

Let $N_{iter}$ represent the iteration count for solving the algebraic system. When the wave number $\omega$ increases (and the mesh size $h$ decreases), the iteration count $N_{iter}$ also increases. In order to describe the growth rate of the iteration count $N_{iter}$ with respect to the wave number $\omega$, we introduce a new notation $\rho$. Let $\omega _1$ and $\omega _2$ be two wave numbers, and let $N^{(1)}_{iter}$ and $N^{(2)}_{iter}$ denote the corresponding iteration counts, respectively. Then we define the positive number $\rho$ by
$$ (\frac {\omega _2}{\omega _1})^\rho=\frac {N^{(2)}_{iter}}{N^{(1)}_{iter}}.$$
For example, when $\rho=1$, the growth is linear; if $\rho\rightarrow 0^+$, then the preconditioner possesses the optimal convergence. For a preconditioner, the positive number $\rho$ defined above is called as ``relative growth rate" of the iteration count.  Of course, we hope that the relative growth rate $\rho$ is sufficiently small. In particular, a preconditioner is almost the optimal if the relative growth rate $\rho$ is
much less than $1$.

\subsection{An example with known analytic solution}

The first model problem is the problem with the Robin boundary condition (refer to \cite{Hut2009}):
\begin{equation}
\begin{split}
&\Delta u+\omega^2u=0 \quad \text{in} \quad \Omega,\\
&{\partial u\over
\partial {\bf n}}+i\omega u=g \quad \text{on}
\quad \partial\Omega,
\end{split}
 \label{eqn1}
\end{equation}
where $\Omega=[0,2]\times[0,1]$, and $g=({\partial \over
\partial {\bf n}}+i\omega)u_{ex}$.

The analytic solution of the problem can be given in the closed form as
$$u_{ex}(x,y)=\text{cos}(k\pi y)(A_1e^{-i\omega_x x}+A_2e^{i\omega_x
x})$$ where $\omega_x=\sqrt{\omega^2-(k\pi)^2}$, and coefficients
$A_1$ and $A_2$ satisfy the equation
\begin{equation}
\left( {\begin{array}{cc} \omega_x & -\omega_x \\
(\omega-\omega_x)e^{-2i\omega_x} & (\omega+\omega_x)e^{2i\omega_x}
\end{array} }
\right)
 \left ( {\begin{array}{c} A_1 \\ A_2
\end{array}}
\right ) = \left ( {\begin{array}{c}
-i   \\
0
\end{array}}
\right )
\end{equation}

In applications, the parameter $k$ may has different values. According to our numerical experiments, different values of $k$ do not affect the efficiency
of the preconditioners (refer to Table 4 and Table 5 in \cite{refhy}). Thus, in order to shorten the length of the paper, we only choose $k=10$ in the experiments for the example.

Let $u_h$ denote the approximate solution generated by an iterative method, we
introduce the following relative error:
$$\text{err.}={||u_{ex}-u_h||_{L^2(\Omega)}\over{||u_{ex}||_{L^2(\Omega)}}} \ .$$
We use the above relative $L^2$ error to measure the accuracy of the
approximate solution $u_h$.

\subsubsection{Results on the PWDG method}

In this part, we apply the PWDG method to the discretzation of this example and solve the resulting algebraic system by PGMRES method, with the preconditioners $B$, $B_s^{(m_0)}$ and $M_1^{(m_0)}$.
In Table \ref{t2}, Table \ref{t10} and Table \ref{t14}, we list the iteration counts and the $L^2$ errors of the resulting approximations.
\vskip 0.1in
\begin{center}
       \tabcaption{
       }  \label{t2}
       PWDG discretization and PGMRES iteration

       (with the~preconditioner $B$)
\vskip 0.1in
\begin{tabular}{|c|c|c|c|c|c|c|} \hline
  $\omega$ & $p$ & $n_h$ & $N_{iter}$ & $\rho$ &  {\text{err.}}  \\ \hline
 $20\pi$ & 10 & $32^2$ & 38 & &  8.13e-4 \\ \hline
 $40\pi$ & 11 & $64^2$ & 47 & 0.3067 &  7.69e-4   \\ \hline
 $80\pi$ & 12 & $128^2$ & 58 & 0.3034 & 6.57e-4   \\ \hline
 $160\pi$ & 15 & $256^2$ & 71 & 0.2918 & 6.02e-4   \\ \hline
 $320\pi$ & 16 & $512^2$ & 87 & 0.2932& 5.98e-4  \\ \hline

 \end{tabular}
 \end{center}

 \vskip 0.1in

\begin{center}
       \tabcaption{
       }  \label{t10}
       PWDG discretization and PGMRES iteration

       (with the~preconditioner $B_s^{(m_0)}$)
\vskip 0.1in
\begin{tabular}{|c|c|c|c|c|c|c|c|c|} \hline
 \multicolumn{3}{|c } {} & \multicolumn{3}{|c}{${m_0=2}$} & \multicolumn{3}{|c|}{${m_0=3}$} \\\hline
  $\omega$ & $p$ & $n_h$ & $N_{iter}$ &$\rho$ &  {\text{err.}} &
 $N_{iter} $  & $\rho$ &{\text{err.}} \\ \hline
 $20\pi$ & 10 & $32^2$ & 44 & & 4.13e-4  &  42 & & 4.23e-4  \\ \hline
 $40\pi$ & 11 & $64^2$ & 53 & 0.2685 & 6.27e-4  &  50 & 0.2515 & 6.21e-4 \\ \hline
 $80\pi$ & 12 & $128^2$ & 64 & 0.2721 & 5.21e-4  &  59 & 0.2388 & 3.97e-4 \\ \hline
 $160\pi$ & 15 & $256^2$ & 77 & 0.2668 & 3.87e-4  &  69 & 0.2259 & 4.27e-4 \\ \hline
 $320\pi$ & 16 & $512^2$ & 92 & 0.2568 & 3.96e-4 & 80 & 0.2134 & 4.12e-4 \\ \hline
 \end{tabular}
 \end{center}

\vskip 0.1in

\begin{center}
       \tabcaption{
       }  \label{t14}
       PWDG discretization and PGMRES iteration

       (with the~preconditioner $M^{(m_0)}_1$)
\vskip 0.1in
\begin{tabular}{|c|c|c|c|c|c|c|c|c|} \hline
 \multicolumn{3}{|c } {} & \multicolumn{3}{|c}{${m_0=2}$} & \multicolumn{3}{|c|}{${m_0=3}$} \\\hline
  $\omega$ & $p$ & $n_h$ & $N_{iter}$ &$\rho$ &  {\text{err.}} &
 $N_{iter} $  & $\rho$ &{\text{err.}} \\ \hline
 $20\pi$ & 10 & $32^2$ & 39 & & 6.52e-4  &  36 & & 6.38e-4  \\ \hline
 $40\pi$ & 11 & $64^2$ & 44 & 0.1740 & 5.87e-4  &  40 & 0.1520 & 5.06e-4 \\ \hline
 $80\pi$ & 12 & $128^2$ & 50 & 0.1844 & 6.28e-4  &  44 & 0.1375 & 5.39e-4 \\ \hline
 $160\pi$ & 15 & $256^2$ & 56 & 0.1635 & 6.14e-4  &  48 & 0.1255 & 5.22e-4 \\ \hline
 $320\pi$ & 16 & $512^2$ & 63 & 0.1699 & 6.29e-4  &  52 & 0.1155 & 7.81e-4 \\ \hline
 \end{tabular}
 \end{center}

\vskip 0.1in

The results in the above tables indicate that the proposed preconditioners are robust for Helmholtz equation with large wave numbers (some detailed comments will be given later).

\subsubsection{Results on the PWLS method}

In this part, we apply the PWLS method to the discretzation of this example and solve the resulting systems by PCG method, with the preconditioners $B$, $B_s^{(m_0)}$, $M_2^{(m_0)}$ and $M_3^{(m_0)}$.
We report the iteration counts and the $L^2$ errors of the resulting approximations in the following four tables.
\vskip 0.1in
\begin{center}
       \tabcaption{
       }  \label{t1}
       PWLS discretization and PCG iteration

       (with the~preconditioner $B$)
\vskip 0.1in
\begin{tabular}{|c|c|c|c|c|c|c|} \hline
  $\omega$ & $p$ & $n_h$ & $N_{iter}$ & $\rho$ & {\text{err.}}  \\ \hline
 $20\pi$ & 10 & $32^2$ & 41 & &  9.25e-4 \\ \hline
 $40\pi$ & 11 & $64^2$ & 51 & 0.3149 & 3.60e-3   \\ \hline
 $80\pi$ & 14 & $128^2$ & 63 & 0.3049 & 3.88e-4   \\ \hline
 $160\pi$ & 15 & $256^2$ & 78 & 0.3081 & 2.31e-4   \\ \hline
 $320\pi$ & 16 & $512^2$ & 96 & 0.2996 & 3.27e-4  \\ \hline

 \end{tabular}
 \end{center}

 \vskip 0.1in

\begin{center}
       \tabcaption{
       }  \label{t9}
       PWLS discretization and PCG iteration

       (with the~preconditioner $B_s^{(m_0)}$)
\vskip 0.1in
 \begin{tabular}{|c|c|c|c|c|c|c|c|c|c|} \hline
 \multicolumn{3}{|c } {} & \multicolumn{3}{|c}{${m_0=2}$} & \multicolumn{3}{|c|}{${m_0=3}$} \\\hline
  $\omega$ & $p$ & $n_h$ & $N_{iter}$ &$\rho$ &  {\text{err.}} &
 $N_{iter} $  & $\rho$ &{\text{err.}} \\ \hline
 $20\pi$ & 10 & $32^2$ & 47 & & 3.64e-4  &  45 & & 3.67e-4  \\ \hline
 $40\pi$ & 11 & $64^2$ & 57 & 0.2783 & 1.79e-3  &  53 & 0.2361 & 1.78e-3 \\ \hline
 $80\pi$ & 14 & $128^2$ & 69 & 0.2756 & 2.63e-4  &  62 & 0.2263 & 3.89e-4 \\ \hline
 $160\pi$ & 15 & $256^2$ & 83 & 0.2665 & 3.91e-4  &  72 & 0.2157 & 2.67e-4 \\ \hline
 $320\pi$ & 16 & $512^2$ & 100 & 0.2688 & 4.37e-4 & 84 & 0.2224 & 4.63e-4 \\ \hline
 \end{tabular}
 \end{center}

 \vskip 0.1in

\begin{center}
       \tabcaption{
       }  \label{t17}
       PWLS discretization and PCG iteration

       (with the~preconditioner $M^{(m_0)}_2$)
\vskip 0.1in
 \begin{tabular}{|c|c|c|c|c|c|c|c|c|} \hline
 \multicolumn{3}{|c } {} & \multicolumn{3}{|c}{${m_0=2}$} & \multicolumn{3}{|c|}{${m_0=3}$} \\\hline
  $\omega$ & $p$ & $n_h$ & $N_{iter}$ &$\rho$ &  {\text{err.}} &
 $N_{iter} $  & $\rho$ &{\text{err.}} \\ \hline
 $20\pi$ & 10 & $32^2$ & 28 & & 6.97e-4  &  26 & & 6.23e-4  \\ \hline
 $40\pi$ & 11 & $64^2$ & 32 & 0.1926 & 2.37e-3  &  29 & 0.1575 & 2.13e-3 \\ \hline
 $80\pi$ & 14 & $128^2$ & 36 & 0.1699 & 4.07e-4  & 32 & 0.1420 & 6.94e-4 \\ \hline
 $160\pi$ & 15 & $256^2$ & 41 & 0.1876 & 6.24e-4  & 35 & 0.1293 & 7.83e-4 \\ \hline
 $320\pi$ & 16 & $512^2$ & 46 & 0.1660 & 4.51e-4  &  38 & 0.1186 & 5.68e-4 \\ \hline
 \end{tabular}
 \end{center}

 \vskip 0.1in

\begin{center}
       \tabcaption{
       }  \label{t13}
       PWLS discretization and PCG iteration

       (with the~preconditioner $M^{(m_0)}_3$)
\vskip 0.1in
 \begin{tabular}{|c|c|c|c|c|c|c|c|c|} \hline
 \multicolumn{3}{|c } {} & \multicolumn{3}{|c}{${m_0=2}$} & \multicolumn{3}{|c|}{${m_0=3}$} \\\hline
  $\omega$ & $p$ & $n_h$ & $N_{iter}$ &$\rho$ &  {\text{err.}} &
 $N_{iter} $  & $\rho$ &{\text{err.}} \\ \hline
 $20\pi$ & 10 & $32^2$ & 23 & & 8.24e-4  &  22 & & 8.29e-4  \\ \hline
 $40\pi$ & 11 & $64^2$ & 26 & 0.1769 & 2.69e-3  &  24 & 0.1255 & 2.70e-3 \\ \hline
 $80\pi$ & 14 & $128^2$ & 29 & 0.1575 & 3.91e-4  &  26 & 0.1155 & 3.90e-4 \\ \hline
 $160\pi$ & 15 & $256^2$ & 33 & 0.1864 & 2.61e-4  &  28 & 0.1069 & 2.34e-4 \\ \hline
 $320\pi$ & 16 & $512^2$ & 37 & 0.1651 & 2.97e-4  &  30 & 0.0995 & 2.28e-4 \\ \hline
 \end{tabular}
 \end{center}

 \vskip 0.1in

It can be seen, from the above tables, that the proposed multilevel preconditioners for Helmholtz equation
with large wave numbers have relatively stable convergence. Namely, the iteration counts of the corresponding iterative methods (PCG or PGMRES) increase slowly when the wave number increases (and the mesh size decreases).
In particular, for the multiplicative multilevel overlapping preconditioners with smoothers, the relative growth rates $\rho$ of the iteration counts with respect to the wave numbers are very small. In fact, the rates are about $0.1$ when the smoothing step $m_0=3$. This means that the multiplicative multilevel overlapping  preconditioners with smoothers are almost optimal. We also notice that, for the PWLS method, the non-standard symmetrized preconditioner $M^{(m_0)}_3$ is more effective than the standardly symmetrized preconditioner $M^{(m_0)}_2$. We would like to emphasize that all the results are obtained without the limiting condition
on the coarsest mesh size $d_0$ (see Section 1 for the details), which can be chosen as a constant independent of $\omega$ and the mesh size $h$.

In the next part, we report some results to explain why the proposed preconditioners are robust for the considered model, and
illustrate the differences between the proposed preconditioners and several existing preconditioners.

\subsubsection{Results on some other related preconditioners}

In this part, we only apply the PWLS method to the discretzation of this example and solve the resulting systems by PCG method with the considered preconditioners.

At first we consider the preconditioners generated by the non-overlapping domain decomposition method, the domain decomposition method with one element overlap and the domain decomposition method with
complete overlap, respectively. Here we consider only the usual one-level domain decomposition (i.e., $J=1$), in which $\Omega$ is decomposed into $4\times 4$ rectangles with the
same size. The resulting preconditioners are denoted by $M_{non}$, $M_{small}$
and $M_{large}$. We give the iteration counts of the PCG methods with the three preconditioners in Table \ref{table-new1}.

\vskip 0.1in

\begin{center}
       \tabcaption{
       }  \label{table-new1}
       PWLS discretization and PCG iteration

       (with the preconditioners $M_{non}$, $M_{small}$ and $M_{large}$)
\vskip 0.1in
\begin{tabular}{|c|c|c|c|c|c|c|c|c|} \hline
 \multicolumn{3}{|c } {} & \multicolumn{2}{|c}{${M_{non}}$} & \multicolumn{2}{|c}{${M_{small}}$}& \multicolumn{2}{|c|}{${M_{large}}$} \\\hline
  $\omega$ & $p$ & $n_h$ & $N_{iter}$ &$\rho$ &
 $N_{iter} $  & $\rho$ & $N_{iter} $  & $\rho$ \\ \hline
 $20\pi$ & 10 & $32^2$ & 108 &   &  82 &  &  20 &   \\ \hline
 $40\pi$ & 11 & $64^2$ & 139 & 0.3641   &  101 & 0.3007  &  21 &0.0704  \\ \hline
 $80\pi$ & 14 & $128^2$ & 178 & 0.3568   &  125 & 0.3076  &  22 & 0.0671 \\ \hline
 $160\pi$ & 15 & $256^2$ & 229 & 0.3635   &  154 & 0.3010  &  23 & 0.0641 \\ \hline
 \end{tabular}
 \end{center}

The above results indicate that, when we decompose $\Omega$ into several subdomains only, all the standard domain decomposition preconditioners have stable convergence (of course, 
the preconditioner with large overlap converges more rapidly).
But, for this one-level decomposition, each subdomain still contains too many fine elements when $h$ is small (i.e., $\omega$ is large). Because of this, we have to design multilevel
domain decomposition in Section 3, such that each considered domain is decomposed into only several subdomains, and every subdomain at the final level contains several fine elements.
Then each local space decomposition (\ref{3.local_decomposition}) is stable, and so the global space decomposition (\ref{decom-new}) should be stable too. This can intuitively explains why the
proposed multilevel preconditioners are effective for Helmholtz equations with large wave numbers.

Then we investigate the influence of the overlapping degree $\theta_0$ to the effectiveness of the multilevel preconditioner defined by (\ref{3.preconditioner}).
When decreasing the thickness of the overlap to be one fine element (i.e., $\theta_0={h\over d}$), the resulting
multilevel preconditioner is denoted by $B_{small}$ (the preconditioner with small overlap). Let $B_{half}$ denote the multilevel preconditioner with $\theta_0={1\over 2}$ (half overlap).
In the table below, we list the iteration counts of the PCG methods with the two preconditioners and the errors of the resulting approximations.
 \vskip 0.1in

\begin{center}
       \tabcaption{
       }  \label{table-new2}
       PWLS discretization and PCG iteration

       (with the~preconditioners $B_{small}$ and $B_{half}$)
\vskip 0.1in
\begin{tabular}{|c|c|c|c|c|c|c|c|c|} \hline
 \multicolumn{3}{|c } {} & \multicolumn{3}{|c}{$B_{small}$} & \multicolumn{3}{|c|}{$B_{half}$} \\\hline
  $\omega$ & $p$ & $n_h$ & $N_{iter}$ &$\rho$ &  {\text{err.}} &
 $N_{iter} $  & $\rho$ &{\text{err.}} \\ \hline
 $20\pi$ & 10 & $32^2$ & 71 & & 3.86e-4  &  49 & & 4.91e-4  \\ \hline
 $40\pi$ & 11 & $64^2$ & 99 & 0.4796 & 5.34e-4  &  61 & 0.3160 & 6.35e-4 \\ \hline
 $80\pi$ & 14 & $128^2$ & 138 & 0.4792 & 2.74e-4  &  76 & 0.3172 & 2.71e-4 \\ \hline
 $160\pi$ & 15 & $256^2$ & 193 & 0.4839 & 1.67e-4  &  95 & 0.3219 & 1.66e-4 \\ \hline

 \end{tabular}
 \end{center}

The above results tell us that the multilevel preconditioner with small overlap is not satisfactory. Fortunately, the multilevel preconditioner with half overlap
possesses almost the same convergence rate as the multilevel preconditioner $B$ with complete overlap (comparing the results in Table \ref{t1}). Notice that the overlap degree
of the small overlap case depends on $h$, but the overlap degree for the case of complete overlap or half overlap is independent of $h$. This means that the convergence rate
of the proposed preconditioner is mainly determined by the overlap degree, as in the standard overlapping domain decomposition method for diffusion equations.

In the following  we compare the proposed preconditioner $B$ with two standard multilevel preconditioners.  Let $\hat{B}$ (MG-Schwarz) be the multilevel preconditioner defined by (\ref{6.preconditioner1}),
and let MG-Jacobi denote the multilevel preconditioner with Jacobi smoothers (see the first part in Subsection 3.3). For the comparison, we use $4\times 4$ refinement for all cases, i.e.,
choosing $N_0=4\times 4$ in Subsection 3.1 and setting $h_j=h_{j-1}/4$ in Subsection 3.3.
We report the iteration counts of the PCG methods with the three preconditioners in Table \ref{table-new3}
 \vskip 0.1in

\begin{center}
       \tabcaption{
       }  \label{table-new3}
       PWLS discretization and PCG iteration

       (with the~preconditioners MG-Jacobi, MG-Schwarz and $B$)
\vskip 0.1in
\begin{tabular}{|c|c|c|c|c|c|c|c|c|c|c|} \hline
 \multicolumn{3}{|c } {} & \multicolumn{2}{|c}{MG-Jacobi} & \multicolumn{2}{|c|}{MG-Schwarz} &\multicolumn{2}{|c|}{${B}$}\\\hline
  $\omega$ & $p$ & $n_h$ & $N_{iter}$ &$\rho$ &
 $N_{iter} $  & $\rho$ & $N_{iter} $  & $\rho$ \\ \hline
 $20\pi$ & 10 & $32^2$ & 78 & &     62 & & 41 &   \\ \hline
 $40\pi$ & 11 & $64^2$ & 113 & 0.5348 &  86 &0.4721 &51 & 0.3124  \\ \hline
 $80\pi$ & 14 & $128^2$ & 163 & 0.5285 &   119 & 0.4686& 63 & 0.3049 \\ \hline
 $160\pi$ & 15 & $256^2$ & 235 & 0.5278 &   164 & 0.4627& 78 & 0.3081 \\ \hline
 \end{tabular}
 \end{center}

\vskip 0.1in
The results given in the above table indicate that the proposed preconditioner $B$ is essentially different from the standard multilevel preconditioners and is obviously more
effective than the considered two preconditioners (see Subsection 3.3 for the detailed explanations). We point out that, when setting $h_j=h_{j-1}/2$ in Subsection 3.3 or implementing
more smoothing steps of the smoothers $B^{-1}_j$ and $\hat{B}^{-1}_j$, this conclusion still holds.

Now we compare three preconditioners, in which each subproblem to be solved has $p$ unknowns. When setting $h_j=h_{j-1}/2$ and implementing $m_0$ smoothing steps for
the Jacobi smoothers, the resulting multigrid preconditioner with Jacobi smoothers is denoted by MG-Jacobi$^{(m_0)}$. If the smoothing step $m_0$ in the preconditioner $B_s^{(m_0)}$
described in Section 4 is not fixed, but it is determined by Krylov method (see \cite{Elam2001}) with the control accuracy $\eta$, the resulting preconditioner is denoted by $B_{s,\eta}$.
As an example, we choose $m_0=3$ and $\eta={1\over 5}$, for which the average time for implementing smoothers in $B_{s,\eta}$ is about $2.7$.
In table \ref{table-new4}, we list the iteration counts of the PCG methods with the three preconditioners.

\vskip 0.1in

\begin{center}
       \tabcaption{
       }  \label{table-new4}
       PWLS discretization and PCG iteration

       (with the~preconditioners MG-Jacobi$^{(m_0)}$, $B_{s,\eta}$ and $B_s^{(m_0)}$,~where $m_0=3,\eta=1/5$)

\vskip 0.1in
\begin{tabular}{|c|c|c|c|c|c|c|c|c|} \hline
 \multicolumn{3}{|c } {} & \multicolumn{2}{|c}{MG-Jacobi$^{(m_0)}$} & \multicolumn{2}{|c}{$B_{s,\eta}$}& \multicolumn{2}{|c|}{$B_s^{(m_0)}$} \\\hline
  $\omega$ & $p$ & $n_h$ & $N_{iter}$ &$\rho$ &
 $N_{iter} $  & $\rho$ & $N_{iter} $  & $\rho$ \\ \hline
 $20\pi$ & 10 & $32^2$ & 61 &   &  48 &  &  45 &   \\ \hline
 $40\pi$ & 11 & $64^2$ & 83 & 0.4443   &  58 & 0.2730  &  53 &0.2361  \\ \hline
 $80\pi$ & 14 & $128^2$ & 112 & 0.4323   &  70 & 0.2713  &  62 & 0.2263 \\ \hline
 $160\pi$ & 15 & $256^2$ & 151 & 0.4310   &  84 & 0.2630  &  72 & 0.2157 \\ \hline
 \end{tabular}
 \end{center}

It can be seen from the above results that the proposed preconditioner $B_s^{(m_0)}$ is obviously more effective than the multigrid preconditioner with $m_0$ Jacobi smoothing steps, and it
is so effective as the preconditioner $B_{s,\eta}$. As pointed out in \cite{Elam2001}, the use of Krylov methods often plays an important role in other methods, but the conclusion
is not true in the current multilevel method.

Notice that we have not reported the errors of the approximations in Table \ref{table-new1}, Table \ref{table-new3} and Table \ref{table-new4} because of the limitation of the space in these tables.
In fact, all the errors are less than $10^{-3}$ and have not large difference.

\subsection{An example whose analytic solution is unknown}

The example tested in the last subsection is too special. In this subsection, we consider the model with
an arbitrary function $g$, which is not determined by an analytic solution. The example can be described as

\begin{equation}
\begin{split}
&\Delta u+\omega^2u=0 \quad \text{in} \quad \Omega,\\
&{\partial u\over
\partial {\bf n}}+i\omega u=g \quad \text{on}
\quad \partial\Omega,
\end{split}
 \label{eqn1}
\end{equation}
where $\Omega=[0,2]\times[0,1]$, and $g=x*y$.

In this example, since we do not know its analytic solution, we can only compute an approximate solution for the comparison with the iterative solution. Let $\hat{u}_h$ be the approximate solution obtained by the direct method for the discrete system, i.e.,
\[
\hat{u}_h=A^{-1}f_h.
\]
To measure the accuracy of the approximate solution $u_h$ generated by an iterative method, we introduce the following relative error:
$$\text{err.}={||\hat{u}_h-u_h||_{L^2(\Omega)}\over{||\hat{u}_h||_{L^2(\Omega)}}} \ .$$

\vskip 0.1in
\subsubsection{Results on the PWDG method}

 In this part we apply the PWDG method to the discretzation of this example and solve the resulting algebraic system by PGMRES method, with the preconditioners $B$, $B_s^{(m_0)}$ and $M_1^{(m_0)}$.
 In Table \ref{t4}, Table \ref{t12} and Table \ref{t16}, we report the iteration counts and the $L^2$ errors of the resulting approximations.
\vskip 0.1in
\begin{center}
       \tabcaption{
       }  \label{t4}
       PWDG discretization and PGMRES iteration

       (with the~preconditioner $B$)
\vskip 0.1in
\begin{tabular}{|c|c|c|c|c|c|c|} \hline
  $\omega$ & $p$ & $n_h$ & $N_{iter}$ &$\rho$ &  {\text{err.}}  \\ \hline
 $20\pi$ & 10 & $32^2$ & 44 & & 7.21e-4 \\ \hline
 $40\pi$ & 11 & $64^2$ & 54 & 0.2955 & 7.39e-4   \\ \hline
 $80\pi$ & 12 & $128^2$ & 66 & 0.2895 & 6.33e-4   \\ \hline
 $160\pi$ & 15 & $256^2$ & 80 & 0.2775 & 6.19e-4   \\ \hline
 $320\pi$ & 16 & $512^2$ & 98 & 0.2928 & 5.63e-4  \\ \hline

 \end{tabular}
 \end{center}

 \vskip 0.1in

\begin{center}
       \tabcaption{
       }  \label{t12}
       PWDG discretization and PGMRES iteration

       (with the~preconditioner $B_s^{(m_0)}$)
\vskip 0.1in
\begin{tabular}{|c|c|c|c|c|c|c|c|c|} \hline
 \multicolumn{3}{|c } {} & \multicolumn{3}{|c}{${m_0=2}$} & \multicolumn{3}{|c|}{${m_0=3}$} \\\hline
  $\omega$ & $p$ & $n_h$ & $N_{iter}$ &$\rho$ &  {\text{err.}} &
 $N_{iter} $  & $\rho$ &{\text{err.}} \\ \hline
 $20\pi$ & 10 & $32^2$ & 51 & & 8.71e-4  &  48 & & 6.22e-4  \\ \hline
 $40\pi$ & 11 & $64^2$ & 61 & 0.2583 & 5.46e-4  &  57 & 0.2479 & 5.83e-4 \\ \hline
 $80\pi$ & 12 & $128^2$ & 73 & 0.2591 & 6.74e-4  &  67 & 0.2332 & 5.91e-4 \\ \hline
 $160\pi$ & 15 & $256^2$ & 87 & 0.2531 & 7.93e-4  & 78 & 0.2193 & 6.08e-4 \\ \hline
 $320\pi$ & 16 & $512^2$ & 104 & 0.2575 & 6.28e-4 & 91 & 0.2224 & 5.69e-4 \\ \hline
 \end{tabular}
 \end{center}

 \vskip 0.1in

\begin{center}
       \tabcaption{
       }  \label{t16}
       PWDG discretization and PGMRES iteration

       (with the~preconditioner $M^{(m_0)}_1$~)
\vskip 0.1in
\begin{tabular}{|c|c|c|c|c|c|c|c|c|} \hline
 \multicolumn{3}{|c } {} & \multicolumn{3}{|c}{${m_0=2}$} & \multicolumn{3}{|c|}{${m_0=3}$} \\\hline
  $\omega$ & $p$ & $n_h$ & $N_{iter}$ &$\rho$ &  {\text{err.}} &
 $N_{iter} $  & $\rho$ &{\text{err.}} \\ \hline
 $20\pi$ & 10 & $32^2$ & 47 & & 7.08e-4  &  41 & & 8.26e-4  \\ \hline
 $40\pi$ & 11 & $64^2$ & 53 & 0.1733 & 6.13e-4  &  45 & 0.1343 & 7.93e-4 \\ \hline
 $80\pi$ & 12 & $128^2$ & 60 & 0.1790 & 6.29e-4  &  49 & 0.1229 & 6.15e-4 \\ \hline
 $160\pi$ & 15 & $256^2$ & 68 & 0.1806 & 6.37e-4  &  53 & 0.1132 & 7.04e-4 \\ \hline
 $320\pi$ & 16 & $512^2$ & 77 & 0.1793 & 7.24e-4  &  57 & 0.1050 & 6.87e-4 \\ \hline
 \end{tabular}
 \end{center}

\vskip 0.1in

The above results indicate that the proposed preconditioners are also robust for this example.

\subsubsection{Results on the PWLS method}

In this part we apply the PWLS method to the discretzation of this example and solve the resulting systems by PCG method, with the preconditioners $B$, $B_s^{(m_0)}$, $M_2^{(m_0)}$ and $M_3^{(m_0)}$.
We list the iteration counts and the $L^2$ errors of the resulting approximations in the following four tables.
\begin{center}
       \tabcaption{
       }  \label{t3}
       PWLS discretization and PCG iteration

       (with the preconditioner $B$)
\vskip 0.1in
\begin{tabular}{|c|c|c|c|c|c|c|} \hline
  $\omega$ & $p$ & $n_h$ & $N_{iter}$ & $\rho$ & {\text{err.}}  \\ \hline
 $20\pi$ & 10 & $32^2$ & 49 & & 6.27e-4 \\ \hline
 $40\pi$ & 11 & $64^2$ & 60 & 0.2922 & 9.19e-3   \\ \hline
 $80\pi$ & 14 & $128^2$ & 73 & 0.2829 & 4.70e-4   \\ \hline
 $160\pi$ & 15 & $256^2$ & 89 & 0.2859 & 3.08e-4   \\ \hline
 $320\pi$ & 16 & $512^2$ & 109 & 0.2925 & 5.81e-4  \\ \hline

 \end{tabular}
 \end{center}
 \vskip 0.1in

\begin{center}
       \tabcaption{
       }  \label{t11}
       PWLS discretization and PCG iteration

       (with the~preconditioner $B_s^{(m_0)}$)
\vskip 0.1in
\begin{tabular}{|c|c|c|c|c|c|c|c|c|} \hline
 \multicolumn{3}{|c } {} & \multicolumn{3}{|c}{${m_0=2}$} & \multicolumn{3}{|c|}{${m_0=3}$} \\\hline
  $\omega$ & $p$ & $n_h$ & $N_{iter}$ &$\rho$ &  {\text{err.}} &
 $N_{iter} $  & $\rho$ &{\text{err.}} \\ \hline
 $20\pi$ & 10 & $32^2$ & 56 & & 7.89e-4  &  53 & & 7.61e-4  \\ \hline
 $40\pi$ & 11 & $64^2$ & 67 & 0.2587 & 2.96e-3  &  63 & 0.2494 & 2.37e-3 \\ \hline
 $80\pi$ & 14 & $128^2$ & 80 & 0.2558 & 3.27e-4  &  74 & 0.2322 & 3.82e-4 \\ \hline
 $160\pi$ & 15 & $256^2$ & 96 & 0.2630 & 5.14e-4  &  86 & 0.2168 & 6.35e-4 \\ \hline
 $320\pi$ & 16 & $512^2$ & 115 & 0.2605 & 6.27e-4 & 100 & 0.2176 & 5.23e-4 \\ \hline
 \end{tabular}
 \end{center}

\vskip 0.1in

\begin{center}
       \tabcaption{
       }  \label{t18}
       PWLS discretization and PCG iteration

       (with the~preconditioner $M^{(m_0)}_2$~)
\vskip 0.1in
\begin{tabular}{|c|c|c|c|c|c|c|c|c|} \hline
 \multicolumn{3}{|c } {} & \multicolumn{3}{|c}{${m_0=2}$} & \multicolumn{3}{|c|}{${m_0=3}$} \\\hline
  $\omega$ & $p$ & $n_h$ & $N_{iter}$ &$\rho$ &  {\text{err.}} &
 $N_{iter} $  & $\rho$ &{\text{err.}} \\ \hline
 $20\pi$ & 10 & $32^2$ & 35 & & 5.81e-4  &  33 & & 6.14e-4  \\ \hline
 $40\pi$ & 11 & $64^2$ & 40 & 0.1926 & 3.19e-3  &  36 & 0.1651 & 3.25e-3 \\ \hline
 $80\pi$ & 14 & $128^2$ & 45 & 0.1699 & 6.34e-4  &  40 & 0.1520 & 5.12e-4 \\ \hline
 $160\pi$ & 15 & $256^2$ & 51 & 0.1806 & 4.88e-4  &  44 & 0.1375 & 3.09e-4 \\ \hline
 $320\pi$ & 16 & $512^2$ & 57 & 0.1605 &  5.33e-4 & 48 & 0.1225 & 5.87e-4 \\ \hline
 \end{tabular}
 \end{center}

\vskip 0.1in

\begin{center}
       \tabcaption{
       }  \label{t15}
       PWLS discretization and PCG iteration

       (with the~preconditioner $M^{(m_0)}_3$~)
\vskip 0.1in
\begin{tabular}{|c|c|c|c|c|c|c|c|c|} \hline
 \multicolumn{3}{|c } {} & \multicolumn{3}{|c}{${m_0=2}$} & \multicolumn{3}{|c|}{${m_0=3}$} \\\hline
  $\omega$ & $p$ & $n_h$ & $N_{iter}$ &$\rho$ &  {\text{err.}} &
 $N_{iter} $  & $\rho$ &{\text{err.}} \\ \hline
 $20\pi$ & 10 & $32^2$ & 31 & & 4.76e-4  &  29 & & 4.91e-4  \\ \hline
 $40\pi$ & 11 & $64^2$ & 35 & 0.1751 & 8.68e-3  &  32 & 0.1420 & 9.01e-3 \\ \hline
 $80\pi$ & 14 & $128^2$ & 39 & 0.1561 & 4.71e-4  &  35 & 0.1293 & 4.70e-4 \\ \hline
 $160\pi$ & 15 & $256^2$ & 44 & 0.1740 & 3.92e-4  &  38 & 0.1186 & 3.91e-4 \\ \hline
 $320\pi$ & 16 & $512^2$ & 49 & 0.1553 & 6.13e-4  &  41 & 0.1096 & 6.19e-4 \\ \hline
 \end{tabular}
 \end{center}

\vskip 0.1in
From the above results, we know that the proposed multilevel preconditioners
are also very effective for the Helmholtz equation considered in this subsection.

\section{Conclusion}
In this paper we have constructed several multilevel preconditioners for the Helmholtz systems generated by the
plane wave discretization (PWLS or PWDG), based on a multilevel overlapping domain decomposition method. In particular,
we have designed multilevel overlapping preconditioners with smoothers, which are almost the optimal. The numerical results
have illustrated that the proposed preconditioners possess nearly stable convergence for the
two-dimensional Helmholtz equations with large wave numbers, without the limiting condition on the coarse mesh size.
In the next work we shall extend the proposed methods (with some modifications) to solving three-dimensional Helmholtz equations with large wave numbers.

\end{document}